\documentclass[a4paper,12pt]{amsart}

\pdfoutput=1
\usepackage{multicol}
\usepackage{multirow}
\usepackage{amsmath,latexsym,amsbsy,amssymb,fancyhdr}
\usepackage{array}
\usepackage{enumerate}
\usepackage{amsthm}
\usepackage{booktabs} 
\usepackage{multirow}
\usepackage{tabularx}
 \usepackage{epsfig} 

\usepackage{psfrag,graphicx}
\usepackage{epstopdf}
\usepackage[latin1]{inputenc}
\usepackage[colorlinks=true]{hyperref}
\hypersetup{urlcolor=blue, citecolor=blue, linkcolor=blue}

\newtheorem{Theorem}{Theorem}[section]

\newtheorem{Example}[Theorem]{Example}

\numberwithin{equation}{section}

\makeatletter
\def\@date{23 December 2015}
\let\insertdate\@date
\makeatother

\newcommand{\ang}[1]{\langle #1 \rangle}
\newcommand{\bb}[1]{\lbrace #1 \rbrace}
\newcommand{\co}{\operatorname{co}}

\newcommand{\inter}{\operatorname{int}}

\newcommand{\rk}{\operatorname{rank}}
\newcommand{\RSU}{\mathcal{R}_\leq}

\newcommand{\R}{{\mathbb R}}

\usepackage{xcolor}
\definecolor{palegreen}{rgb}{0.2,0.6,0.2}

\author{Robert Baier}
\address[Robert Baier]{Mathematisches Institut,
Universit\"at Bayreuth,  
95440 Bayreuth, Germany}
\email{robert.baier@uni-bayreuth.de}

\author{Thuy T. T. Le}
\address[Thuy T.T. Le]{Universit\`a di Padova, Dipartimento di Matematica, via Trieste 63, 35121 Padova, Italy}
\email{lethienthuy@gmail.com}
\thanks{The second author is supported by a PhD fellowship for foreign students at the Universit\`a di Padova funded by Fondazione CARIPARO. This paper was developed while the second author was visiting the Department of Mathematics of the University of Bayreuth}
\begin{document}

\title{Construction of the Minimum Time Function \newline
       Via Reachable Sets of Linear Control Systems. \newline
       Part~2: Numerical Computations}



\date{\today}
\begin{abstract}
 In the first part of this paper we introduced an algorithm that uses
 reachable set approximation to approximate the minimum time function of linear control problems.
 To illustrate the error estimates and to demonstrate differences to other numerical approaches
 we provide a collection of numerical examples which either allow higher order of convergence 
with respect to time discretization or 
 where the continuity of the minimum time function 
 cannot be sufficiently granted, i.e.~we study cases in which the minimum time
 function is H\"older continuous or even discontinuous.
 
\end{abstract}
\maketitle

\markright{\hfill Construction of the  minimum time function via reachable sets Part 2\hfill}
\markleft{\hfill R. Baier, T. T. T. Le \hfill}

\section{Introduction}

This part should serve as a collection of academic test examples for calculating
the minimum time function for several, mainly linear control problems
which were previously discussed in the literature.
The numerical approximation of the minimum time function is based on the calculation of reachable sets 
with set-valued quadrature methods and is described in full details in the first part \cite{BLp1}. 
Many links to other numerical approaches for approximating reachable sets are also given there.

In several examples, we compare the error of the 
minimum time function studying the influence of its regularity and of the smoothness of the 
support functions of corresponding set-valued integrands. We verify the obtained error
estimates which involve time and space discretization.
We will not repeat the importance, applications and study of the minimum
time function in this second part and refer also for the notations to the first part~\cite{BLp1}.

We first consider several linear examples with various target 
and control sets and study different levels of regularity of the corresponding
minimum time function. 
The control sets are either one- or two-dimensional polytopes (a segment
or a square) or balls and are varied to study different regularity
allowing high or low order of convergence for the underlying set-valued quadrature method.
In all linear examples, we apply a set-valued combination method of order 1 and 2 (the set-valued 
Riemann sum combined with Euler's method resp.~the set-valued trapezoidal rule with Heun's method). 
For the nonlinear example in Subsection~\ref{subsec_nonlin}, we would like to approximate the 
time-reversed dynamics 
of \eqref{example5} 
 directly by Euler's and Heun's method. This example demonstrates  
that this approach is not restricted to 
 the  class of linear control systems. 
 Although first numerical experiences are gathered, its theoretical justification has to be
gained by a forthcoming paper. 
  
In Subsection~\ref{subsec_non_exp_prop} 
one example demonstrates the need of the strict expanding property of (union of) reachable
sets for characterizing boundary points of the reachable set via time-minimal points (compare~\cite[Propositions~2.18 and~2.19]{BLp1}).
The two-dimensional system in Example~\ref{ex:counter_ex_1}
is not normal and has only a one-dimensional reachable set. 

The section ends with a collection of examples which either are more challenging for numerical 
calculations or partially violate~\cite[Assumptions~2.13 (iv) and~(iv)' in Proposition 2.19]{BLp1}.
Finally, a discussion of our approach and possible improvements can be found in Section~\ref{sec:concl}.

\section{Numerical tests}
\label{sec:num_tests}

The following examples in the next subsections are devoted to illustrating the performance of the error behavior of our proposed approach. 

The space discretization follows the presented approach in 
\cite[Subsec.~3.2]{BLp1} 
and uses supporting points in directions
 \begin{align*}
    l^k & := \bigg( \cos\bigg(2 \pi \frac{k-1}{N_{\mathcal{R}}-1}\bigg), \ \sin\bigg(2 \pi \frac{k-1}{N_{\mathcal{R}}-1}\bigg) \bigg)^\top, \ k=1,\ldots,N_{\mathcal{R}} \\
    \eta^r & := \begin{cases}
                   -1 + 2(r-1) & \quad \text{if $U=[-1,1]$},\ r=1,\ldots,N_U, \\
                   l^r  & \quad \text{if $U \subset \R^2$},\  r=1,\ldots,N_U \\
                \end{cases}
 \end{align*}
and normally choose either $N_U = 2$ for one-dimensional control sets or $N_U = N_{\mathcal{R}}$ for $U \subset \R^2$  in the discretizations of the unit sphere~\cite[(3.11)]{BLp1}. 

The comparison of the two applied methods is done by computing the error with respect to the  
$L^{\infty}$-norm  of the difference between the approximate and the true minimum time 
function evaluated at test points. 
The true minimum time function is delivered 
 analytically by tools  from control theory. 
The  test grid points are distributed uniformly over the domain 
$\mathcal{G}=[-1,1]^2$ with step size $\Delta x= 0.02$.

\subsection{Linear examples}
In the linear, two-dimensional, time-invariant Examples~\ref{ex:1}--\ref{ex:3b} we can
check~\cite[Assumption~2.13~(iv)]{BLp1} 
 \begin{quote}
    $\mathcal{R}(t)$ is \emph{strictly expanding} on the compact interval $[t_0,t_f]$, i.e.~$\mathcal{R}(t_1) \subset \inter \mathcal{R}(t_2)$ for all $t_0\le t_1<t_2\le t_f$
 \end{quote}
in several ways. From the numerical calculations
we can observe this property in the shown figures for the fully discrete reachable sets.
Secondly, we can use the available analytical formula for the minimum time function
resp.~the reachable sets or check 
the Kalman rank condition
 \begin{align*}
    \rk\Big[ B, A B \Big] = 2
 \end{align*}
for time-invariant systems if the target is the origin (see~\cite[Theorems~17.2 and~17.3]{HL}).

We start with an example having a Lipschitz continuous minimum time function and verify the
error estimate in~\cite[Theorem~3.7]{BLp1}. Observe that the numerical error here is only
contributed by the spatial discretization of the target set or control set.

\begin{Example} 
\label{ex:1}
Consider the control dynamics , see \cite{BFS,GL},
\begin{equation}\label{example1}
\dot{x}_1=u_1,\,\,\dot{x}_2=u_2,\,\, (u_1,u_2)^\top \in U \text{ with $U:= B_1(0)$ or $U := [-1,1]^2$ }.
\end{equation}
We consider either the small ball $B_{0.25}(0)$ or the origin as target set $\mathcal{S}$. This is a simple time-invariant example with $\bar A=\begin{bmatrix}
     0&  0  \\[0.3em]
      0 & 0 
\end{bmatrix}$, $\bar B=\begin{bmatrix}
     -1&  0  \\[0.3em]
      0 & -1 
\end{bmatrix}$.
Its fundamental solution matrix is the identity matrix, therefore
\begin{equation*}
\mathcal{R}(t)=\Phi(t,t_0)S+\int_{t_0}^{t}\Phi(t,s)\bar B(s)U=S+(t_0-t)U,
\end{equation*}
and any method from (I)--(III) gives the exact solution, i.e.
$$
  \mathcal{R}_h(t)=\mathcal{R}(t) 
  =S+(t-t_0)U
$$ 
due to the symmetry of $U$. For instance, the set-valued Euler scheme with $h=\frac{t_{j+1}-t_j}{N}$ yields 
\begin{equation*}
\begin{cases}
\mathcal{R}_h(t_{j+1})=\mathcal{R}_h(t_j)+h(\bar A \mathcal{R}_h(t_j)+\bar B U)=\mathcal{R}_h(t_j)-hU,\\
\mathcal{R}_h(t_0)=S,
\end{cases}
\end{equation*}
therefore, $\mathcal{R}_{h}(t_N)=S-NhU=S+( t_N -t_0)U$ and the error is only due to the space discretizations $\mathcal{S}_\Delta \approx \mathcal{S}$, $U_\Delta \approx U$ and does not depend on $h$ 
(see Table~\ref{tab:1}). The error would be the same for finer step size $h$ and $\Delta t$
in time or if a higher-order method is applied. Note that the error for the origin as target
set (no space discretization error) is in the magnitude of the rounding errors of floating
point numbers.


We choose $t_f = 1,\,K=10$ and $N=2$ for the computations. The set-valued 
Riemann sum combined with Euler's method is used.



It is easy to check that the minimum time function is Lipschitz continuous, since 
one of the equivalent Petrov conditions in~\cite{P}, \cite[Chap.~IV, Theorem~1.12]{BCD} with $U=B_1(0)$ or $[-1,1]^2$ hold:
 \begin{align*}
    0 & > \min_{(u_1,u_2)^{\top}\in U} \ang{\nabla d(x,\mathcal{S}),(u_1,u_2)^\top}, \\
    0 & \in \inter\bigg( \bigcup_{u \in U} f(0,u) \bigg) \quad\text{with $f(x,u) = A x + B u$.}
 \end{align*}

Moreover, the support function with respect to the time-reversed dynamics \eqref{example1} 
 \begin{align*}
    \delta^* (l,\Phi(t,\tau)\bar B(\tau)U) & = \begin{cases}
                                                  \|l\| & \quad\text{if $U = B_1(0)$}, \\
                                                  |l_1|+|l_2| & \quad\text{if $U = [-1,1]^2$}
                                               \end{cases}
 \end{align*}
is constant  with respect to the time $t$, so it is trivially arbitrarily continuously differentiable with respect to $t$ with bounded derivatives  uniformly for all $l\in S_{n-1}$. 



\begin{table}[h]
\begin{tabular}{|c|c|c|c|}
	\hline
          $  N_{\mathcal{R}} = N_U $ & $U=B_1(0)$,
             & $U=[-1,1]^2$,
             & $U=[-1,1]^2$, \rule{0ex}{3ex} \\
        &   $\mathcal{S}=B_{0.25}(0)$ 
             & $\mathcal{S}=B_{0.25}(0)$  
             & $\mathcal{S}=\bb{0}$ \rule[-1.5ex]{0ex}{1.5ex} \\
	\hline
  $100$& $ 6.14\times 10^{-4}$  & $ 4.9 \times 10^{-4} $  
     & $ 8.9 \times 10^{-16} $ \rule{0ex}{3ex} \\
     	\hline
   $50$ & $ 24\times 10^{-4}$  & $ 19 \times 10^{-4} $  
      & $ 8.9 \times 10^{-16} $ \rule{0ex}{3ex} \\
      \hline
     $25$ &  $ 0.0258 $  & $ 0.0073  $  
      & $ 8.9 \times 10^{-16} $ \rule{0ex}{3ex} \\
	\hline
\end{tabular}\\[2ex]
\caption{error estimates for Example~\ref{ex:1} with different control and target sets} 
\label{tab:1}
\end{table}
\vspace{1ex}

In Fig.~\ref{fig:1_ball} the minimum time functions are plotted 
for Example~\ref{ex:1} for two different control sets $U = B_1(0)$ (left) and $U = [-1,1]^2$ (right) with the same two-dimensional target set $\mathcal{S} = B_{0.25}(0)$. 
The minimum time function is in general not differentiable everywhere. Since it is 
zero in the interior of the target, one has at most Lipschitz continuity at
the boundary of $\mathcal{S}$.
In Fig.~\ref{fig:1_square_target_origin} the minimum time function is plotted 
for the same control set as in Fig.~\ref{fig:1_ball} (right), but this time the target set
is the origin and not a small ball.

\begin{figure}[htp]
\begin{center}
\includegraphics[scale=0.475]{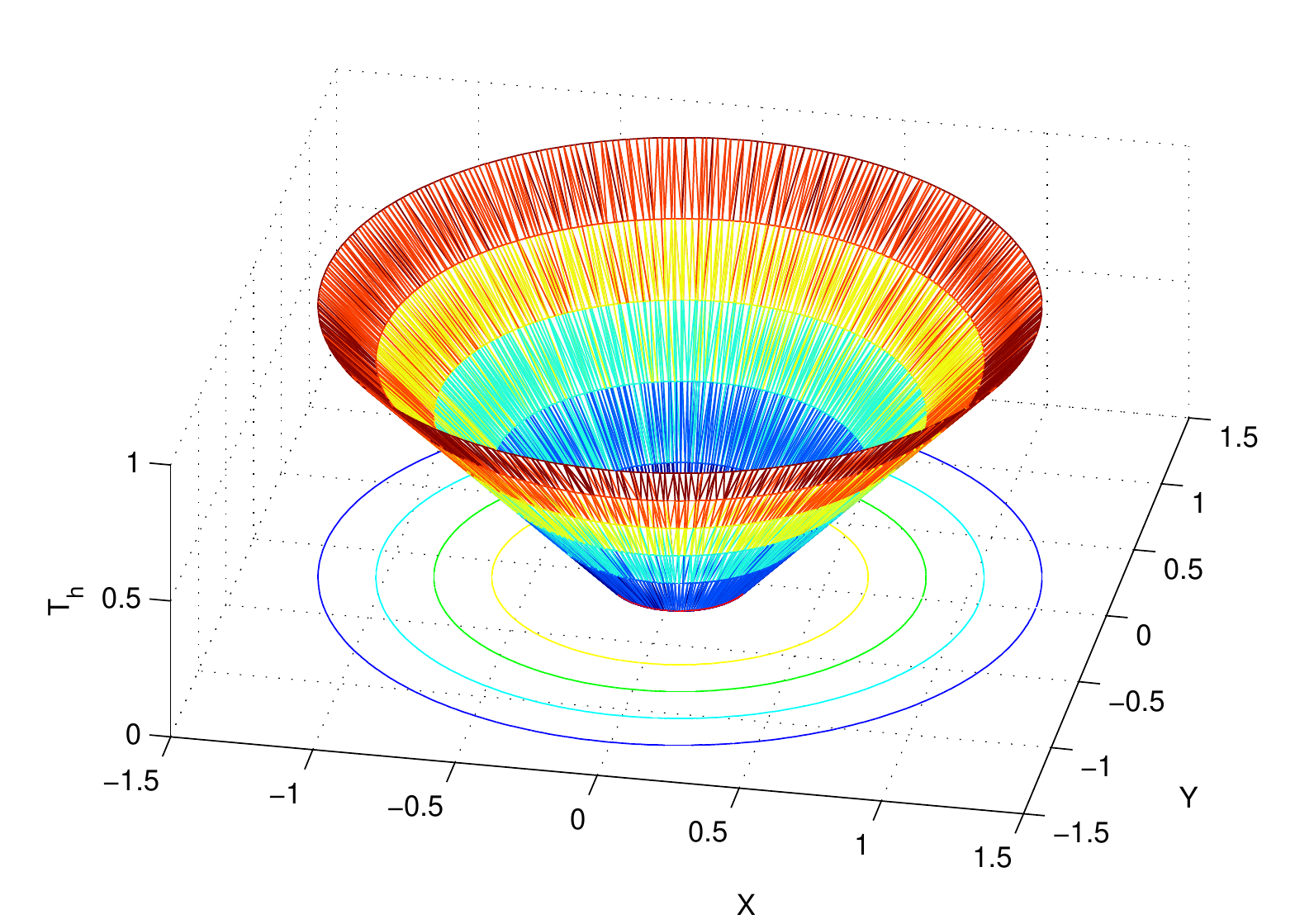}
\includegraphics[scale=0.475]{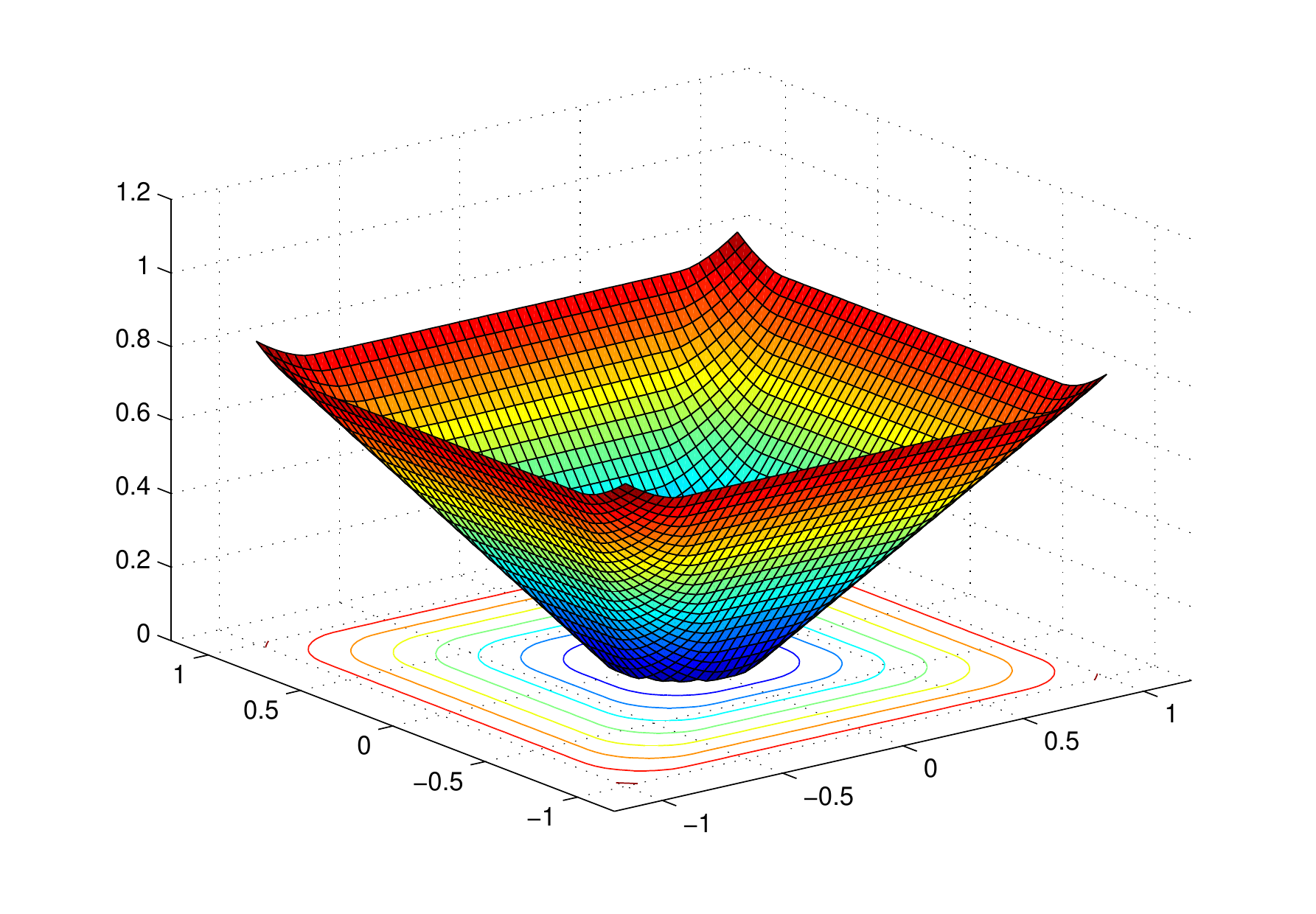}
\caption{minimum time functions for Example~\ref{ex:1} with different control sets}
\label{fig:1_ball}
\end{center}
\end{figure}

\begin{figure}[htp]
\begin{center}
\includegraphics[scale=0.475]{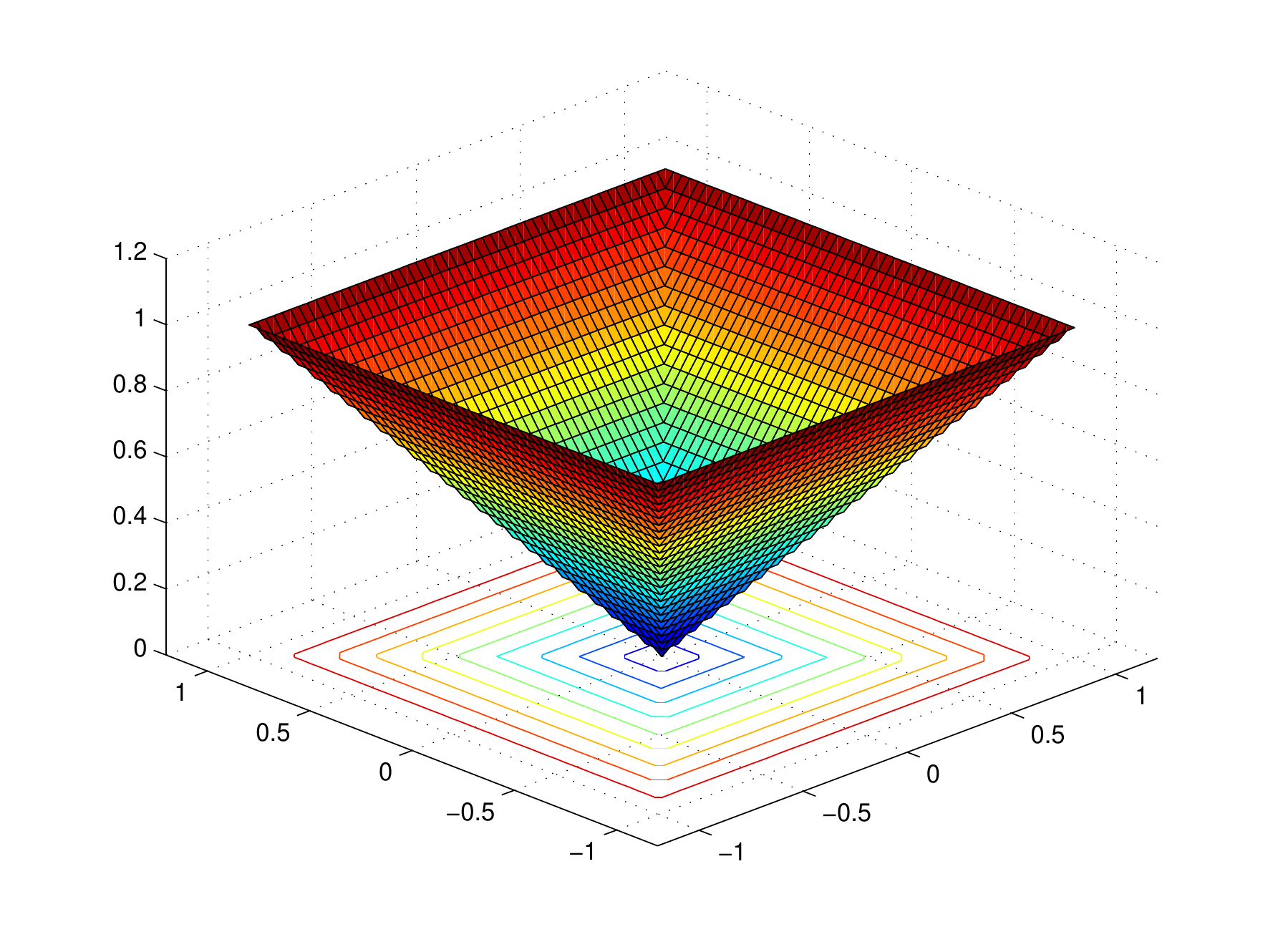}
\caption{minimum time function for Example~\ref{ex:1} with $U = [-1,1]^2$, $\mathcal{S}=\{0\}$}
\label{fig:1_square_target_origin}
\end{center}
\end{figure}
\end{Example} 

We now study well-known dynamics as the double integrator and the harmonic oscillator
in which the control set is one-dimensional. The classical rocket car example with
H\"older-continuous minimum time function was already computed
by the Hamilton-Jacobi-Bellman approach in~\cite[Test~1]{F} and \cite{CL,GL}, where numerical calculations
are carried out by enlarging the target (the origin) by a small ball.

\begin{Example} 
\label{ex:2}
a) The following dynamics is the \emph{double integrator}, see e.g.~\cite{CL}.
\begin{equation}\label{example3}
\dot{x}_1=x_2,\,\dot{x}_2=u,\,\,u\in U := [-1,1].
\end{equation}
We consider either the small ball $B_{0.05}(0)$ or the origin as target set $\mathcal{S}$.
Then the minimum time function is  $\frac{1}{2}$--H\"older continuous 
for the first choice of $\mathcal{S}$ see \cite{AM,CL} and the support function for the time-reversed dynamics \eqref{example3} 
$$\delta^* (l,\Phi(t,\tau)\bar B(\tau)[-1,1])=\delta \Bigg(l,\begin{bmatrix}
     1&  -(t-\tau)  \\[0.3em]
      0 &  1 
\end{bmatrix} \begin{bmatrix}
       0  \\[0.3em]
        -1 
\end{bmatrix}[-1,1]\Bigg)=\big|(t-\tau,-1) \cdot l \big|$$
is only absolutely continuous with respect to $\tau$ for some directions $l \in S_1$
with $l_1 \neq 0$. Hence, we can expect that the convergence order 
for the set-valued quadrature method is at most $2$. \\
We fix $t_f = 1$ as maximal computed value for the minimum time function
and $N = 5$.

 In Table~\ref{tab:2} the error estimates for two set-valued combination methods
are compared (order 1 versus order 2). Since the minimum time function is only
$\frac{1}{2}$--H\"older continuous we expect as overall convergence order $\frac{1}{2}$
resp.~$1$. A least squares approximation of the function $C h^{p}$ for the error term
reveals $C = 1.37606$, $p = 0.4940$ for Euler scheme combined with set-valued Riemann sum
resp.~$C = 22.18877$, $p = 1.4633$ (if $p=1$ is fixed, then $C= 2.62796$) for Heun's method combined with set-valued trapezoidal
rule. Hence, the approximated error term is close to the expected one 
by~\cite[Theorem~3.7 and Remark~3.8]{BLp1}. 
Very similar results are obtained with the Runge-Kutta methods of order 1 and 2
in Table~\ref{tab:22} in which the set-valued Euler method is slightly better than the
combination method of order 1 in Table~\ref{tab:2}, and the set-valued Heun's method
coincides with the combination method of order 2, since both methods use the same approximations of the given dymanics. 

Here we have chosen to double the number of directions $N_{\mathcal{R}}$ each time the step size
is halfened which is suitable for a first order method. For a second order method
we should have multiplied $N_{\mathcal{R}}$ by 4 instead. From this point it is not surprising
that there is no improvement of the error in the fifth row for step size $h = 0.0025$. 

\begin{table}[h]
\begin{tabular}{|l|c|c|c|}
	\hline
	\mbox{ }\ \,$h$ & $N_{\mathcal{R}}$ & \textbf{Euler scheme \& Riemann sum} 
          & \textbf{Heun's scheme \& trapezoid rule} \\
      \hline
        $0.04$ & $50$ & $0.2951$ & $0.2265$ \\
	  \hline
        $0.02$ & $100$ & $0.1862$ & $0.1180$ \\
	  \hline
        $0.01$ & $200$ & $0.1332$ & $0.0122$ \\
	  \hline
        $0.005$ & $400$ & $0.1132$ & $0.0062$ \\
	  \hline
        $0.0025$ & $800$ & $0.0683$ & $0.0062$ \\
	  \hline
\end{tabular}\\[2ex]
\caption{error estimates for Ex.~\ref{ex:2} a) 
         for combination methods of order 1 and 2} 
\label{tab:2}
\end{table}
\begin{table}[h]
\begin{tabular}{|l|c|c|c|}
	\hline
	\mbox{ }\ \,$h$ & $N_{\mathcal{R}}$ & \textbf{set-valued Euler method} 
          & \textbf{set-valued Heun method} \\
          \hline
          	  $0.04$ & $50$ & $0.2330$ & $0.2265$ \\
	      \hline
	          $0.02$ & $100$ & $0.1681$  &  $0.1180$ \\
	  	  \hline
	          $0.01$ & $200$ & $0.1149$   & $0.0122$  \\
	  	  \hline
	          $0.005$ & $400$ & $0.0753$ &  $0.0062$ \\
	  	  \hline
	          $0.0025$ & $800$ &$0.0318$  &  $0.0062$ \\
	  	  \hline
\end{tabular}\\[2ex]
\caption{error estimates for Ex.~\ref{ex:2} a) 
         for Runge-Kutta meth.\ of order 1 and 2} 
\label{tab:22}
\end{table}
%
	
\vspace{1ex}
\begin{figure}[htp]
\begin{center}
\includegraphics[scale=0.5]{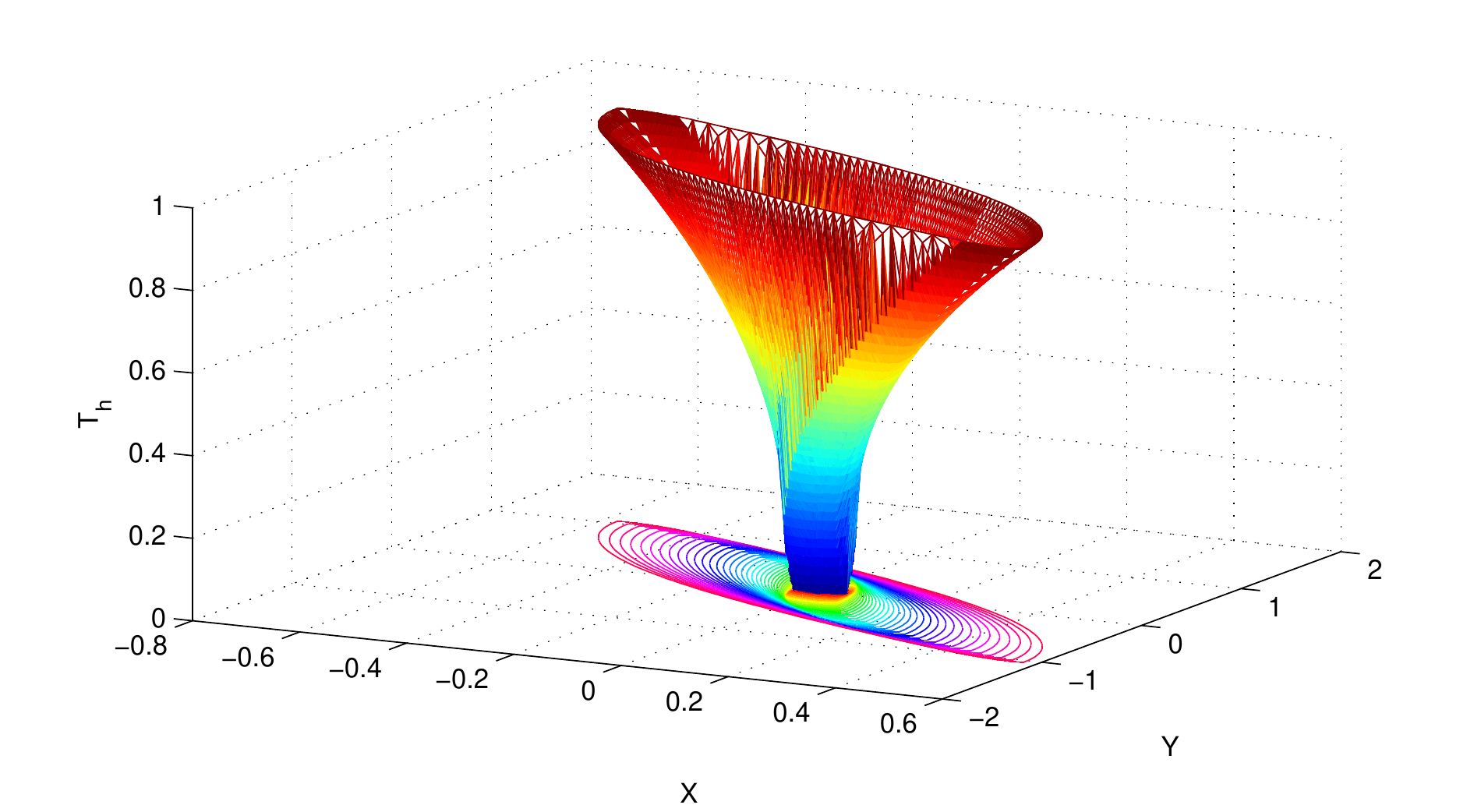}
\caption{minimum time function  for Example~\ref{ex:2}a)  with target set $B_{0.05}(0)$}\label{exam_31}
\label{fig:2a_small_ball}
\end{center}
\end{figure}
 As in Example~\ref{ex:1}  we can consider the dynamics \eqref{example3} with the origin as a target (see the minimum time function in Fig.~\ref{fig:2a_origin}~(left). In this case, the numerical computation by PDE approaches, i.e.~the solution of the associated Hamilton-Jacobi-Bellman equation (see e.g.~\cite{F})  requires the replacement of the target point $0$ by a small ball  $B_\varepsilon(0)$ for  suitable  $\varepsilon>0$. This replacement surely increases the error of the  calculation   (compare the minimum time function in Fig.~\ref{fig:2a_small_ball} for $\varepsilon = 0.05$). However, the proposed approach works perfectly regardless of the fact whether  $\mathcal{S}$  is a  two-dimensional set or a singleton.
\\[1ex]

b) harmonic oscillator dynamics (see~\cite[Chap.~1, Section~1.1, Example 3]{LM})
\begin{equation}\label{example4}
\dot{x}_1=x_2,\,\dot{x}_2=-x_1+u,\,\,u\in U :=[-1,1].
\end{equation} 
Since the Kalman rank condition
 \begin{align*}
    \rk\Big[ B, A B \Big] = 2,
 \end{align*}
the minimum time function $T(\cdot)$ is also continuous.
The plot for $T(x)$ for the harmonic oscillator with the origin as
target, $ t_f=6,\,N_{\mathcal{R}} = 100,\,
N=5 $ and $K=40$ is shown in Fig.~\ref{fig:2a_origin}~(right).

According to \cite[Sec.~3.4]{BLp1} we construct open-loop time-optimal controls for the discrete 
problem with target set $\mathcal{S} = \{0\}$ by Euler's method. In Fig.~\ref{fig:exam_32} the corresponding discrete open-loop time-optimal
trajectories for Examples~\ref{ex:2}a)~(left) and b)~(right) are depicted.
%
%

\begin{figure}[htp]
\begin{center}
\includegraphics[scale=0.5]{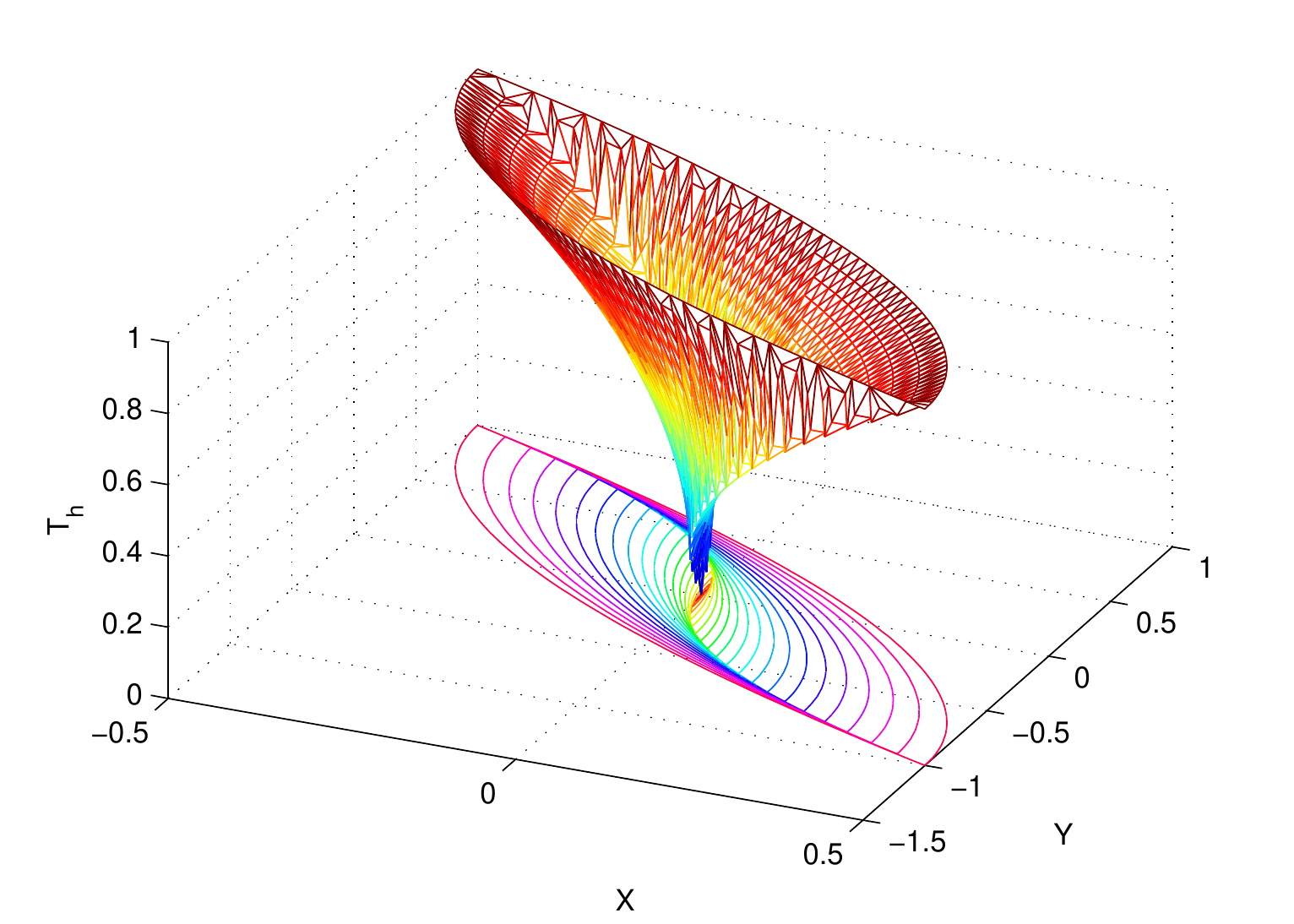}
\ \,
\includegraphics[scale=0.35]{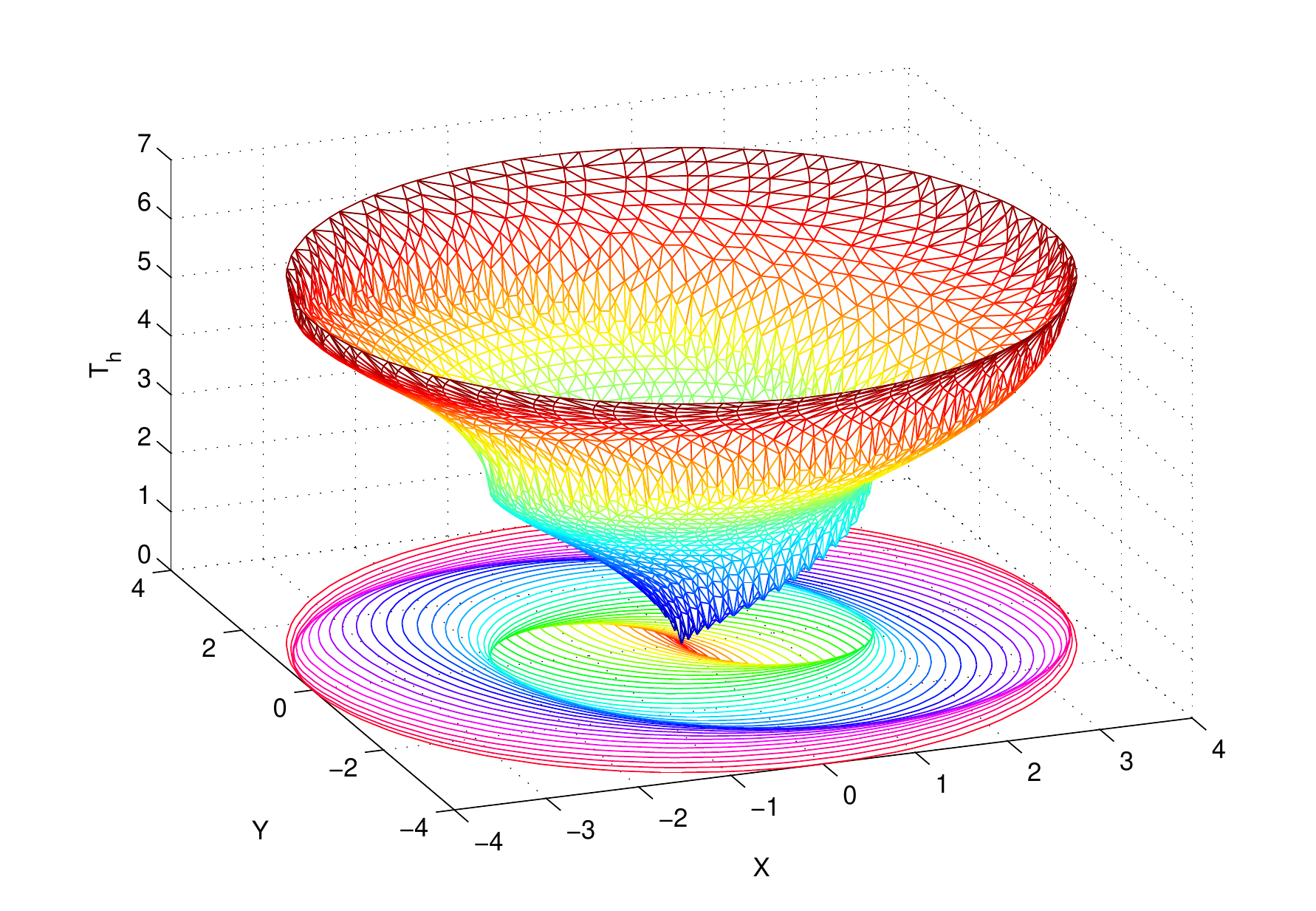}
\caption{minimum time functions for Example~\ref{ex:2}a) resp.~b)}
\label{fig:2a_origin}
\end{center}
\end{figure}


\begin{figure}[htp]
\begin{center}
\includegraphics[width=2.4in,height=2.4in]{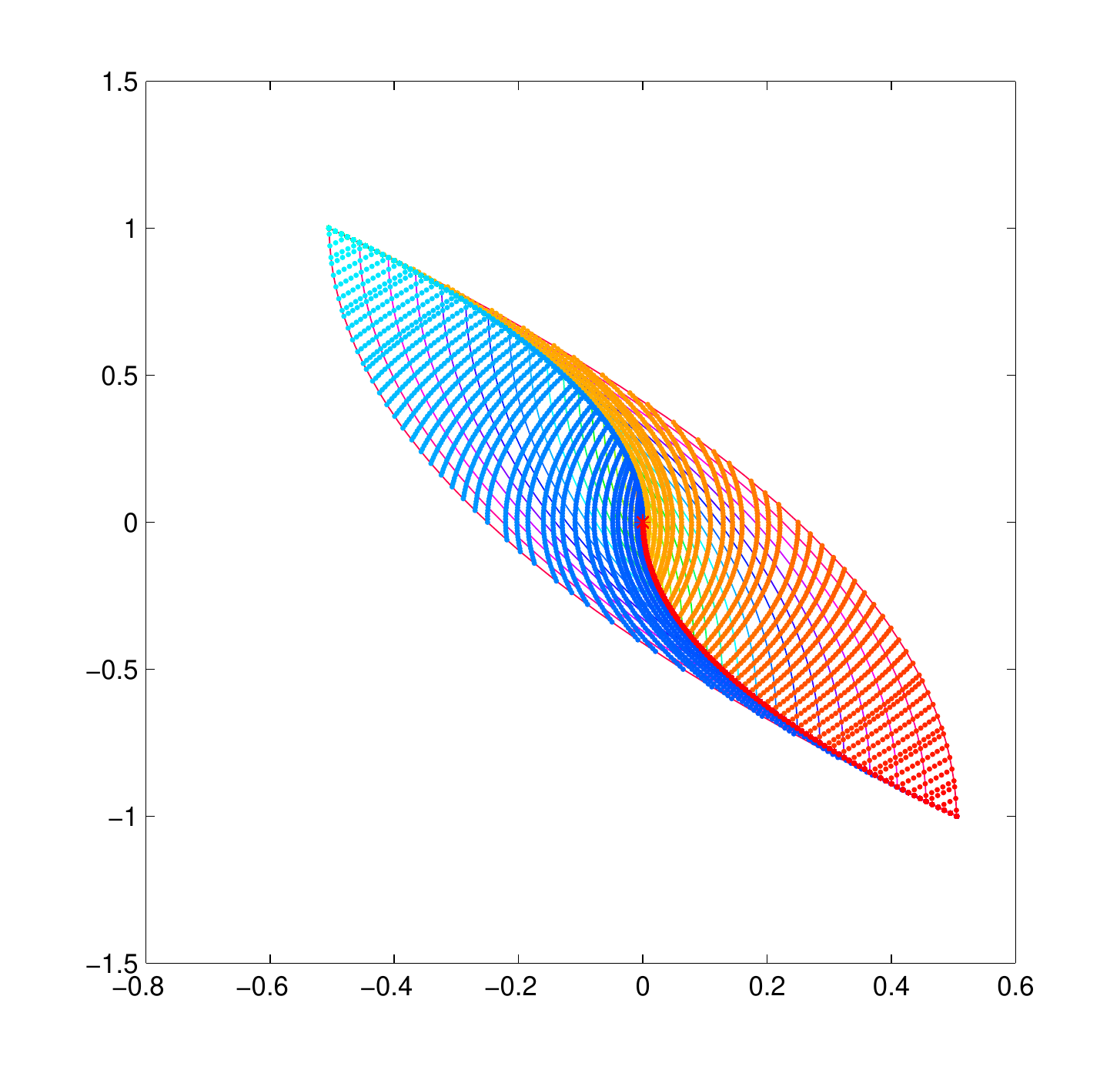}
\includegraphics[width=3.125in,height=2.4in]{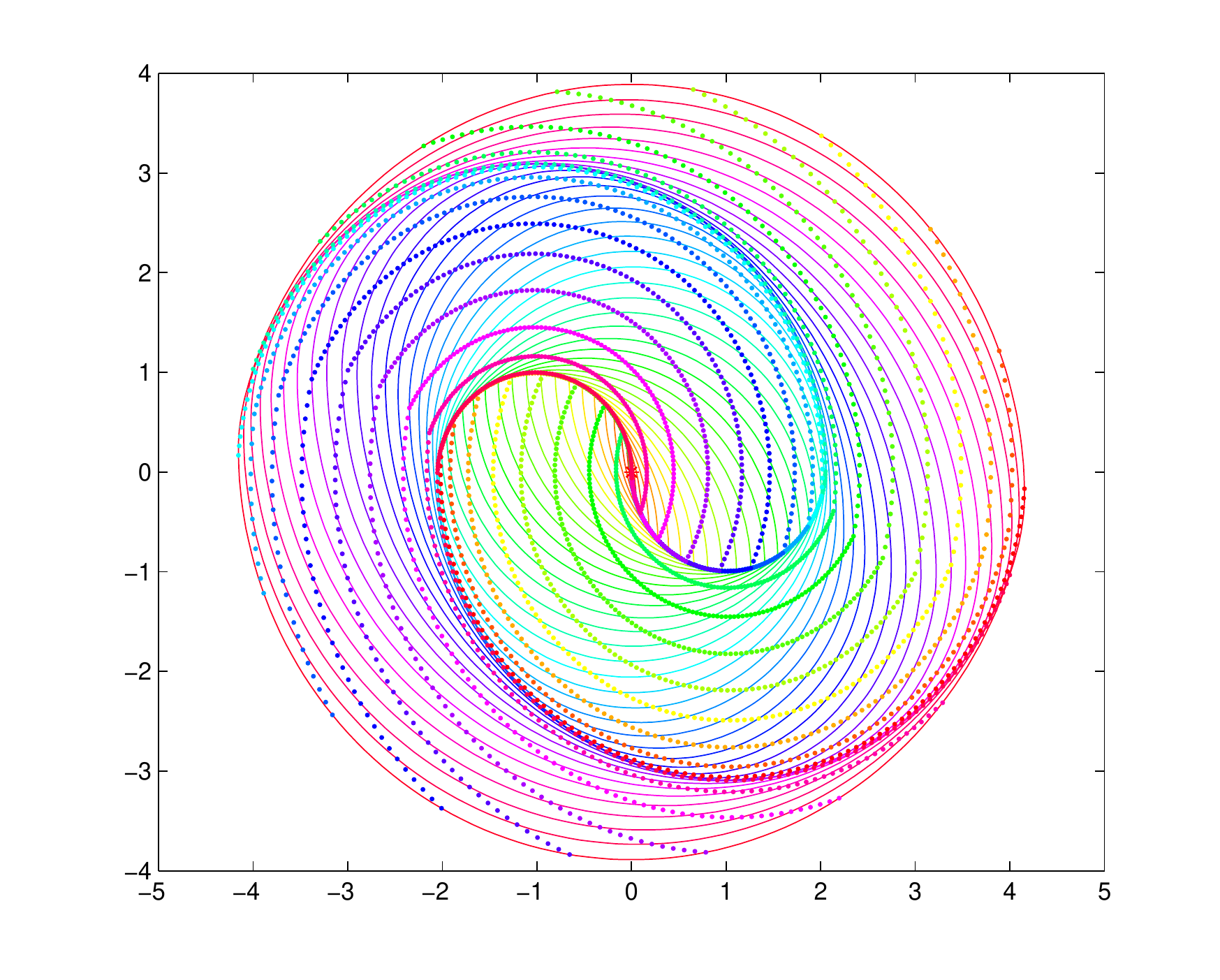}
\caption{approximate optimal trajectories for Example~\ref{ex:2}a) resp.~b)}\label{fig:exam_32}
\end{center}
\end{figure}
\end{Example} 

The following two examples exhibit smoothness of the support functions and would even allow 
for methods with order higher than two with respect to time discretization.
The first example has a special linear dynamics and is smooth, although the control set is a unit square.

\begin{Example} 
\label{ex:3a} 
   In the third linear two-dimensional example the reachable set for various end times $t$ 
   is always a polytope with four vertices and coinciding outer normals at its faces. 
   Therefore, it is a smooth example which would even justify the use of methods with higher order than 2 to
   compute the reachable sets (see \cite{BLcham,BL}).
   It is similar to Example~\ref{ex:counter_ex_1}, but has an additional column in matrix $B$ and is a variant of \cite[Example~2]{BLcham}.

   Again, we fix $t_f = 1$ as maximal time value and compute the result with $N=2$.
   We choose $N_{\mathcal{R}} = 50$ normed directions, since the reachable set has only four different
   vertices.


\begin{equation}\label{example15}
\begin{bmatrix}
    \dot x_1 \\[0.3em]
     \dot x_2 
\end{bmatrix}=\begin{bmatrix}
    0 & -1 \\[0.3em]
    2 &  3  
\end{bmatrix}\begin{bmatrix}
    x_1 \\[0.3em]
      x_2  
\end{bmatrix}+\begin{bmatrix}
    1 & -1 \\[0.3em]
    -1 & 2  
\end{bmatrix}\begin{bmatrix}
    u_1 \\[0.3em]
      u_2  
\end{bmatrix},
\end{equation}
where $(u_1,\,u_2)^\top \in [-1,1]^2$. Let the origin be the target set  $\mathcal{S}$. 
The fundamental solution matrix of the time-reversed dynamics of \eqref{example15} is given by
 \begin{align*}
    \Phi(t,\tau) & = \begin{bmatrix}
                        2 e^{-(t-\tau)} - e^{-2(t-\tau)} 
                          & e^{-(t-\tau)} - e^{-2(t-\tau)} \\[0.3em]
                        -2 e^{-(t-\tau)} + 2 e^{-2(t-\tau)} 
                          & -e^{-(t-\tau)} + 2 e^{-2(t-\tau)}
                     \end{bmatrix}.
 \end{align*}
This is a smooth example in the sense that the support function  for  the time-reversed
set-valued dynamics of \eqref{example15}, 
 \begin{align*}
    \delta ^*(l,\Phi(t,\tau)\bar B(\tau)[-1,1]^2)=e^{-(t-\tau)}\vert l_1-l_2 \vert + e^{-2(t - \tau)}\vert l_1-2l_2 \vert,
 \end{align*}
is smooth with respect to $\tau$  uniformly in $l \in S_1$ .

The analytical formula for the (time-continuous) minimum time function is as follows:
\begin{equation*}
\begin{aligned}
T((x_1,x_2)^\top)=\max \bb{& t \colon \ t \geq 0 \text{ is the solution of one of the equations }\\
 & x_2=-2x_1\pm (e^{-t}-1),\,x_2=-x_1\pm 1/2(1-e^{-2t})}
\end{aligned}
\end{equation*}

A least squares approximation of the function $C h^{p}$ for the error term
reveals $C = 2.14475$, $p = 0.8395$ for the set-valued combination method of order~1
and $C = 23.9210$, $p= 1.7335$ (if $p=2$ is fixed, then $C=70.1265$)
 for the one of order~2. The values are similar to the expected ones from \cite[Remark~3.8]{BLp1}, 
since the minimum time function (see Fig.~\ref{fig:3a}~(left)) is Lipschitz (see~\cite[Sec.~IV.1, Theorem~1.9]{BCD}).

Similarly, another variant of this example with a one-dimensional control can be constructed by deleting
the second column in matrix $B$. The resulting (discrete and continuous-time) reachable sets
would be line segments. Thus, the algorithm would compute the fully discrete minimum time
function on this one-dimensional subspace. The absence of interior points in the reachable
sets is not problematic for this approach in contrary to common approaches based on the
Hamilton-Jacobi-Bellman equation as shown in Example~\ref{ex:counter_ex_1}.
\begin{table}[h]
\begin{tabular}{|l|c|c|}
	\hline
	\mbox{ }\ \ \,$h$ & \textbf{Euler scheme \& Riemann sum} 
          & \textbf{Heun's scheme \& trapezoid rule} \\
	  \hline
        $0.05$ & $0.170\phantom{0}$ & $0.1153\phantom{000}$ \\
	  \hline
        $0.025$ & $0.095\phantom{0}$ & $0.0470\phantom{000}$ \\
	  \hline
        $0.0125$ & $0.0599$ & $0.0133\phantom{000}$ \\
	  \hline
        $0.00625$ & $0.0285$ & $0.0032\phantom{000}$ \\
	  \hline
\end{tabular}\\[2ex]
\caption{error  estimates for Example~\ref{ex:3a} 
         for methods of order 1 and 2} 
\label{tab:3}
\end{table}
\vspace{1ex}

\begin{figure}[htp]
\begin{center}
\includegraphics[scale=0.45]{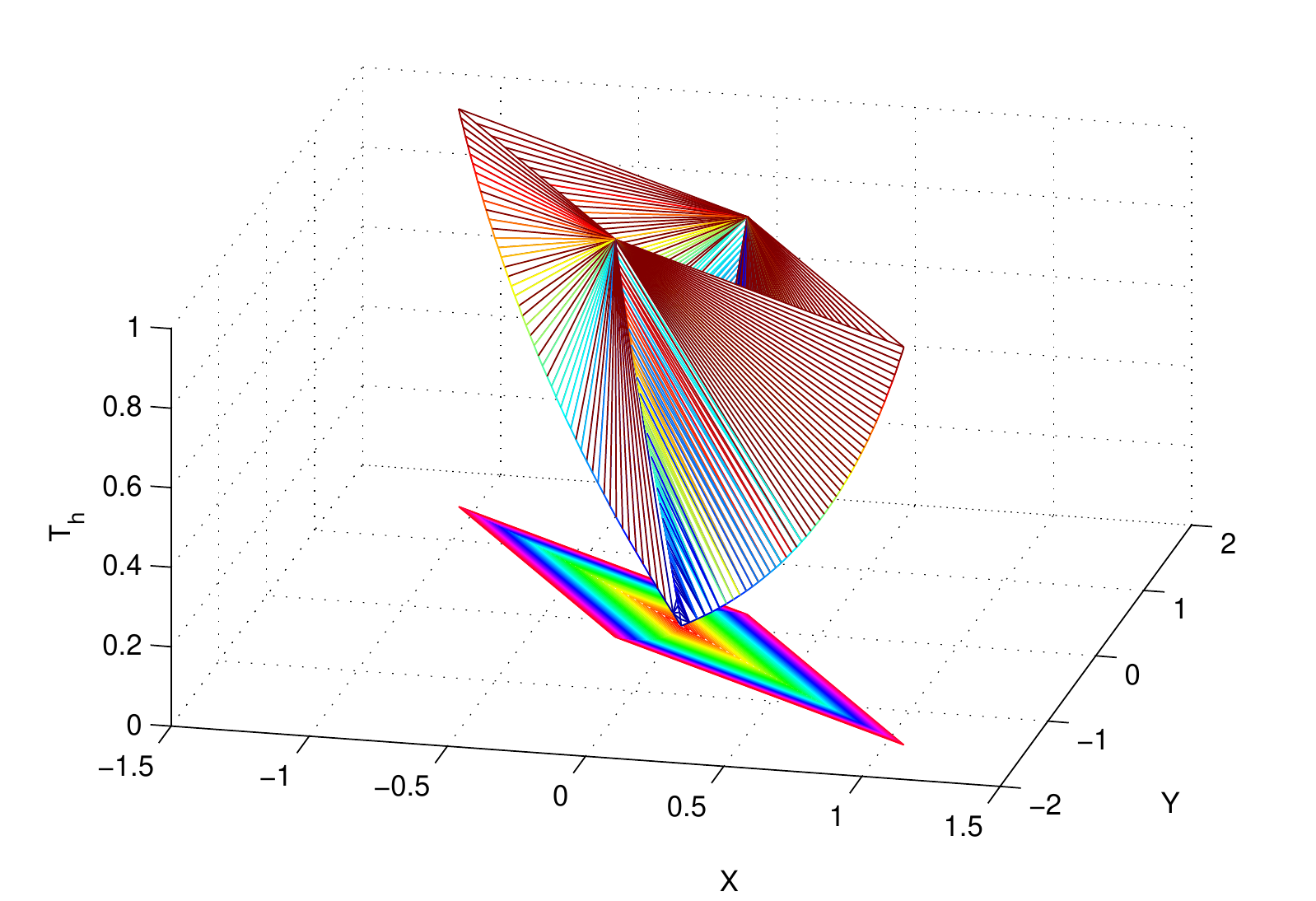}
\ \,
\includegraphics[scale=0.5]{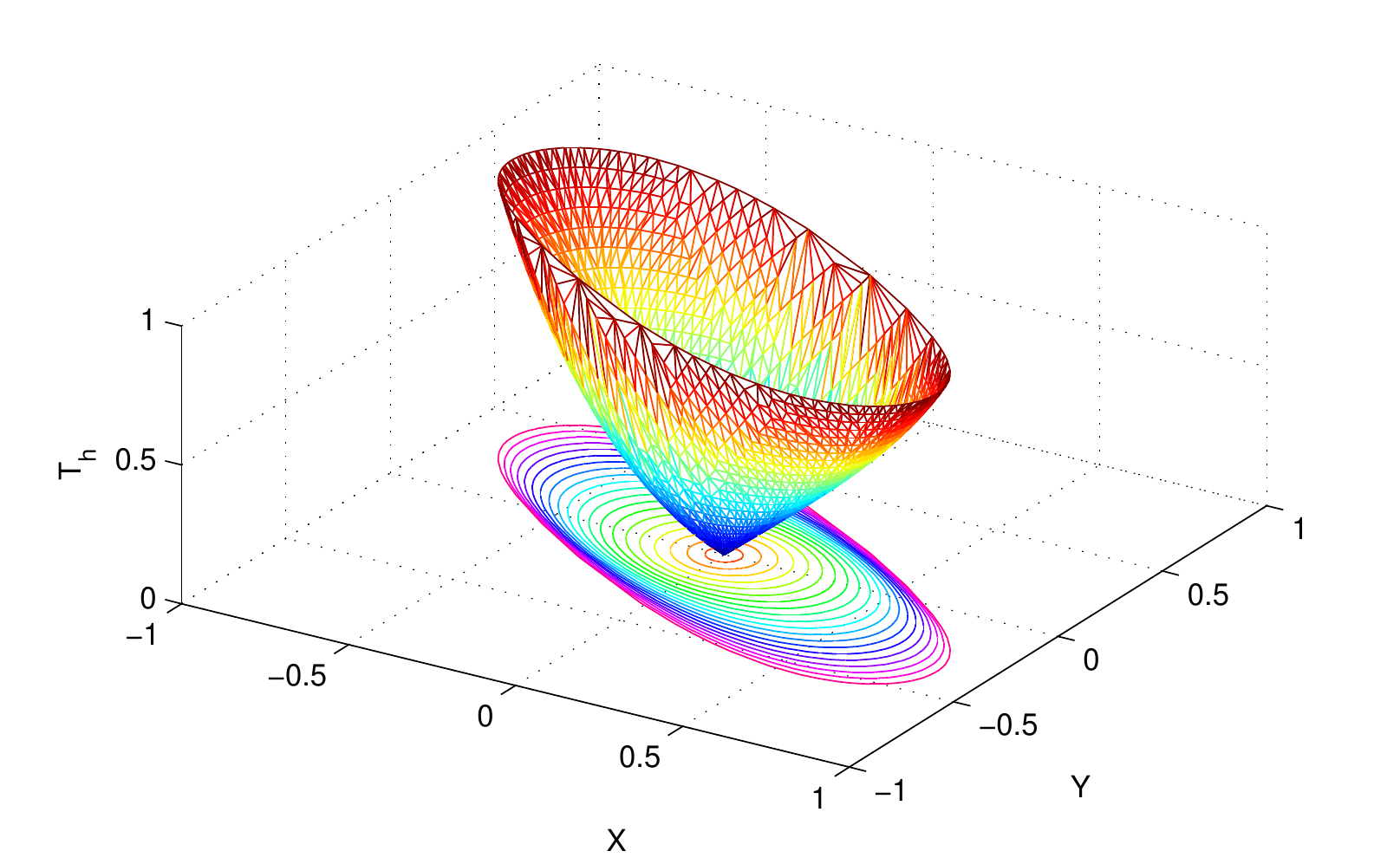}
\caption{ minimum time functions for Examples~\ref{ex:3a} and~\ref{ex:3b}} 
\label{fig:3a}
\end{center}
\end{figure}
\end{Example} 

The next example involves a ball as control set and leads naturally to a smooth support function. 

\begin{Example} 
\label{ex:3b}
  The following smooth example is very similar to the previous example.
It  is given in \cite[Example~4.2]{BBCG}, \cite[Example~4.4]{BL} 
\begin{equation}\label{example16}
\begin{bmatrix}
    \dot{x}_1 \\[0.3em]
    \dot{x}_2  
\end{bmatrix}=\begin{bmatrix}
    0 & -1 \\[0.3em]
    2 & 3  
\end{bmatrix}\begin{bmatrix}
    x_1 \\[0.3em]
      x_2  
\end{bmatrix}+B_1(0)
\end{equation}
and uses a ball as control set. This is a less academic example than Example~\ref{ex:3a} (in which the matrix $B(t)$ was carefully chosen), since a ball as control
set often allows the use of higher order methods for the computation of reachable sets
(see~\cite{BL,B}). Here, no analytic formula for the minimum time function is available
so that we can study only numerically the minimum time function (see Fig.~\ref{fig:3a}~(right)).
Obviously, the support function is again smooth with respect to $\tau$ uniformly in all
normed directions $l$, since
 \begin{align*}
    \delta ^*(l,\Phi(t,\tau) B_1(0) & = \| \Phi(t,\tau)^\top l \|.
 \end{align*}

\begin{figure}[htp]
\begin{center}
\end{center}
\end{figure}
\end{Example} 

\subsection{A nonlinear example}
\label{subsec_nonlin}

The following special bilinear example with convex reachable sets may provide the hope to extend 
our approach to some class of nonlinear dynamics.

\begin{Example} 
\label{ex:4}
The nonlinear dynamics is one of the examples in \cite{GL}.
\begin{equation}\label{example5}
\dot{x}_1=-x_2+x_1 u,\,\,\dot{x}_2=x_1+x_2 u,\,\,u\in [-1,1].
\end{equation}
With this dynamics, after computing the true minimum time function we observe that $T(\cdot)$ is Lipschitz continuous and its sublevel set, which is exactly the reachable set at the corresponding time, satisfies the required properties. The target set  $\mathcal{S}$  is $B_{0.25}(0)$.  \\

We fix $t_f = 1$ as maximal computed value for the minimum time function
and $N=2$.
Estimating the error term $C h^{p}$ in Table~\ref{tab:4} by least squares approximation yields the values
$C = 0.3293133$, $p= 1.8091$ for the set-valued Euler method
and $C = 0.5815318$, $p= 1.9117$ for the Heun method. 

The unexpected good behavior of Euler's method stems from the specific behavior of trajectories.
Although the distance of the end point of the Euler iterates for halfened step size to the
true end point is halfened, but the distance of the Euler iterates to the boundary of the 
true reachable set is almost shrinking by the factor 4 due to the specific tangential
approximation. 
In Fig.~\ref{fig:4} the Euler iterates are marked with \textasteriskcentered\ in red color,
while Heun's iterates are shown with $\circ$ marks in blue color. The symbol $ \bullet $ marks the end point of the corresponding true solution.

Observe that the dynamics originates from the following system in polar coordinates
\begin{equation*} 
\dot{r}=r u,\,\,\dot{\varphi}=1 ,\,\,u\in [-1,1].
\end{equation*}
Hence, the reachable set will grow exponentially with increasing time.

 \begin{figure}[htp]
 \begin{center}
 \includegraphics[scale=0.4]{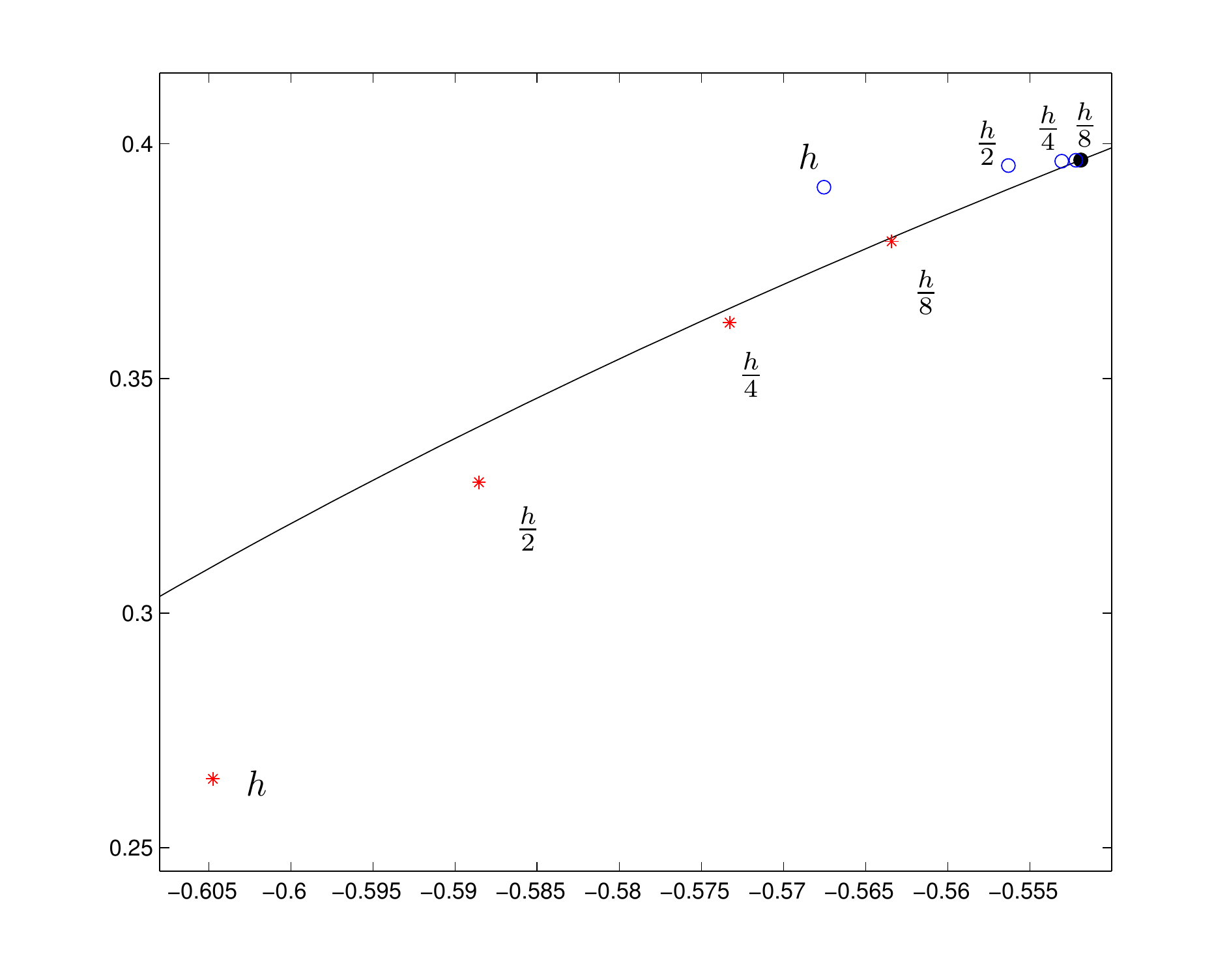}
 \caption{ Euler and Heun's iterates for Example~\ref{ex:4}}  
 \label{fig:4}
 \end{center}
 \end{figure}
  
\begin{table}[h]
\begin{tabular}{|l|c|c|c|}
	\hline
	\mbox{ }\ \,$h$ & $N_{\mathcal{R}}$  & \textbf{set-valued Euler scheme} 
          & \textbf{set-valued Heun's scheme} \\
	  \hline
	$0.5$ & $ 50 $&$0.0848\phantom{00}$ & $0.1461\phantom{00}$ \\
	  \hline
	$0.1$ & $100 $&$0.0060\phantom{00}$ & $0.0076\phantom{00}$ \\
	  \hline
        $0.05$ &$200$  & $0.0015\phantom{00}$ & $0.0020\phantom{00}$ \\
	  \hline
        $0.025$& $400$ & $0.00042\phantom{0}$ & $0.000502$ \\
	  \hline
        $0.0125$ &$800$ & $0.000108$ & $0.000126$ \\
	  \hline
\end{tabular}\\[2ex]
\caption{error estimates for Example~\ref{ex:4} with set-valued methods
         of order 1 and 2}
\label{tab:4}
\end{table}
\vspace{1ex}

The minimum time function for this example is shown in Fig.~\ref{exam_5}.
\begin{figure}[htp]
\begin{center}
\includegraphics[width=4.5in,height=3in]{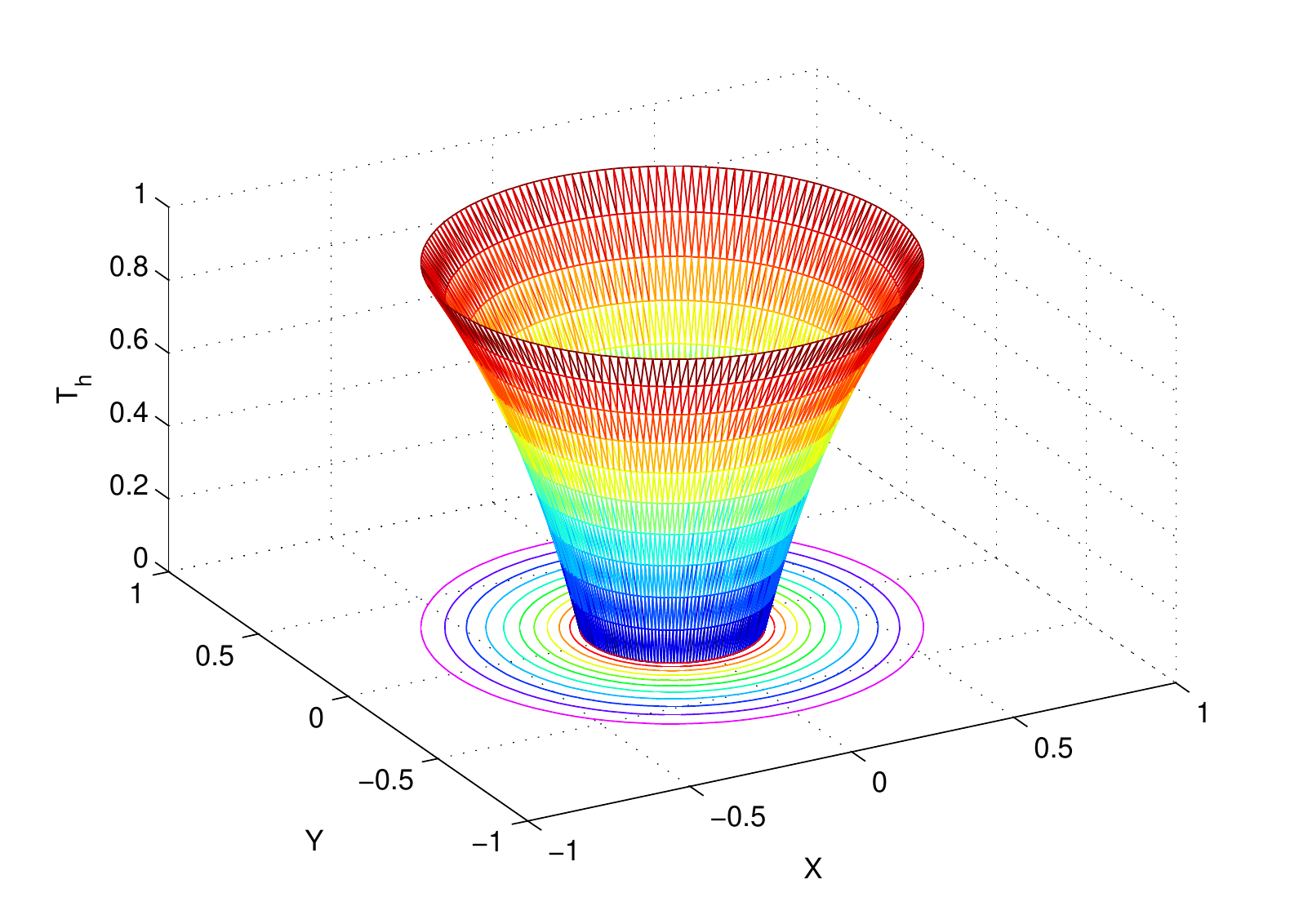}
\caption{minimum time function for Example~\ref{ex:4}} 
\label{exam_5}
\end{center}
\end{figure}
\end{Example} 

\subsection{Non-strict expanding property of reachable sets}
\label{subsec_non_exp_prop}%

The next example violates the continuity of the minimum time function (the dynamics is not normal).
Nevertheless, the proposed~\cite[Algorithm~3.4]{BLp1} is able to provide a good
approximation of the discontinuous minimum time function.

\begin{Example} 
   Consider the dynamics
\label{ex:counter_ex_1}
\begin{equation}
\begin{bmatrix}
    \dot{x}_1 \\[0.3em]
    \dot{x}_2  
\end{bmatrix}=\begin{bmatrix}
    0 & -1 \\[0.3em]
    2 & 3  
\end{bmatrix}\begin{bmatrix}
    x_1 \\[0.3em]
      x_2  
\end{bmatrix}+ u_1 \begin{bmatrix}
     1  \\[0.3em]
    -1   
\end{bmatrix}
\end{equation}
with $u_1 \in U = [-1,1]$, $\mathcal{S} = \{0\}$ and $t \in I=[0, t_f]$.

The Kalman rank condition yields
 \begin{align*}
    \rk\Big[ B, A B \Big] = 1 < 2
 \end{align*}
so that the normality of the system is not fulfilled.

The fundamental system $\Phi(t,\tau)$ (for the time-reversed system)
is the same as in Example~\ref{ex:3a} so that 
 \begin{align*}
    & \delta ^*(l,\Phi(t,\tau)\bar B(\tau)[-1,1])=e^{\tau-t}\vert l_1-l_2 \vert 
      = e^{\tau-t} \delta ^*(l, V), \\
    \intertext{with the line segment $V = \co(\begin{bmatrix}
                                        -1  \\[0.3em]
                                         1   
                                      \end{bmatrix}, \begin{bmatrix}
                                         1  \\[0.3em]
                                        -1   
                                      \end{bmatrix})$. 
                                      Since}
    & \int_0^t \delta^*(l, \Phi(t,\tau)\bar B(\tau)[-1,1]) d\tau 
      = e^{\tau-t} \bigg|_{\tau=0}^t \cdot \delta ^*(l, V) \\
    = & (1 - e^{-t}) \cdot \delta ^*(l, V) = \delta ^*(l,  (1 - e^{-t}) V), \\
    \mathcal{R}(t) & = \int_0^t \Phi(t,\tau)\bar B(\tau)[-1,1] d\tau =  (1 - e^{-t}) V.
 \end{align*}
Hence, the reachable set is an increasing line segment (and always part of the same line
in $\R^2$, i.e.~it is one-dimensional so that the interior is empty). Clearly,
both inclusions 
 \begin{align} \label{ex:relaxed_expand}
    \mathcal{R}(s) \subset \mathcal{R}(t) \quad\text{or}\quad
    \RSU(s) \subset \RSU(t)
 \end{align}
i.e.~(see~\cite[(2.12)]{BLp1}),  
hold, but not the strictly expanding property of $\overline{\mathcal{R}}(\cdot)$ 
on $[t_0, t_f]$ in~\cite[Assumptions~2.13~(iv) and~(iv)']{BLp1}, i.e.
 \begin{align}
    \overline{\mathcal{R}}(t_1) & \subset \inter \overline{\mathcal{R}}(t_2) \text{\ for all $t_0\le t_1<t_2\le t_f$, where \rule{0ex}{4ex}} \label{eq:strict_exp_prop} \\
    %
        \overline{\mathcal{R}}(t) & = \begin{cases}
                                         \mathcal{R}(t) & \text{for Assumption (iv)}, \\
                                         \RSU(t) & \text{for Assumption (iv)'}.
                                      \end{cases}\nonumber
     \rule[-5.5ex]{0ex}{5.5ex}
 \end{align}
By \cite[Sec.~IV.1, Proposition~1.2]{BCD} the minimum time function is discontinuous (it has infinite values outside the line segment).

The plots of the two continuous-time reachable sets $\mathcal{R}(t)$ for $t=1,2$ together with the true minimum time function (in red)
and its discrete analogue (in green) obtained by the Euler scheme with $h = 0.025$ are shown
in Fig.~\ref{fig:counter_exam}:
\begin{figure}[htp]
\begin{center}
\includegraphics[scale=0.4]{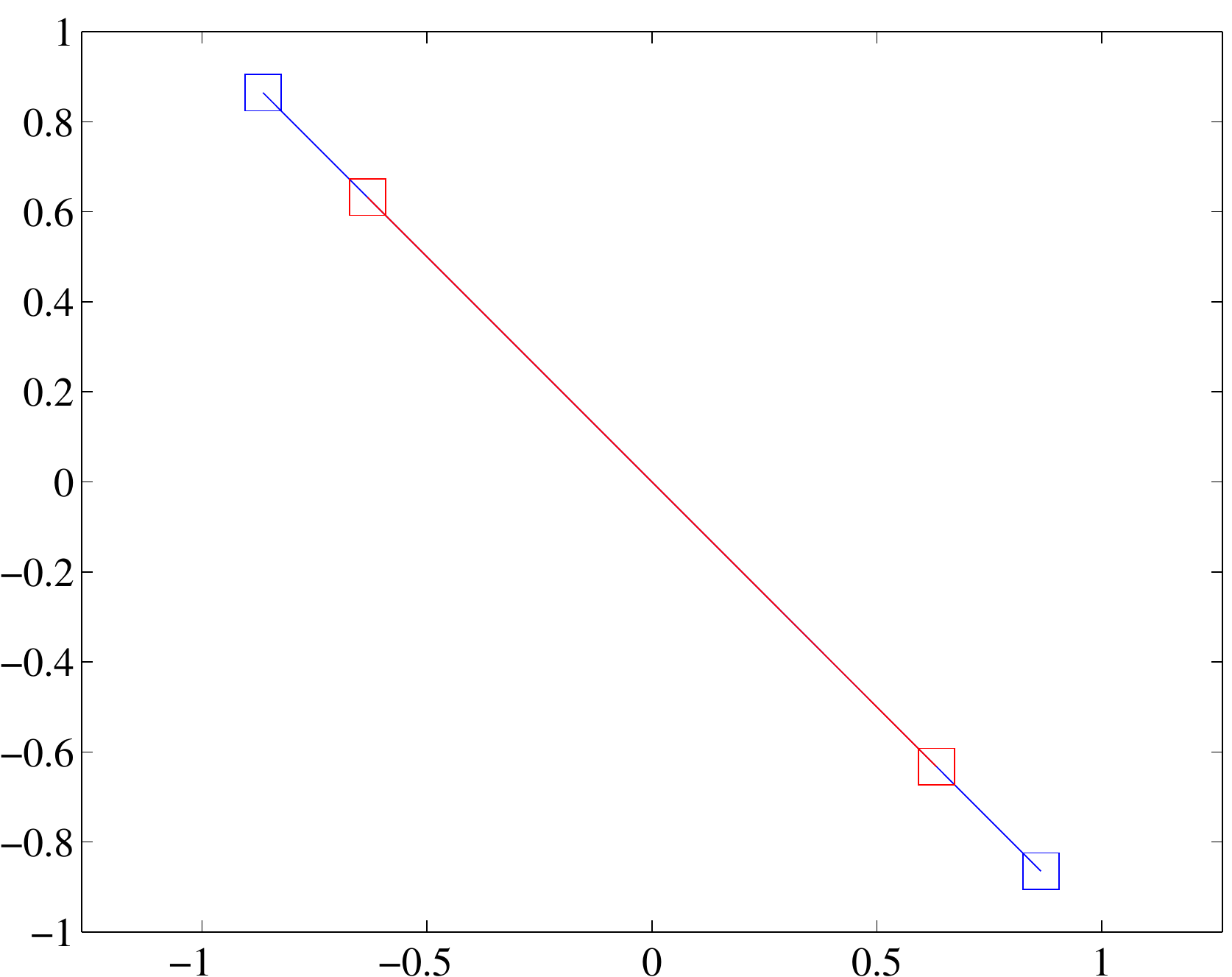}
\quad 
\includegraphics[scale=0.45]{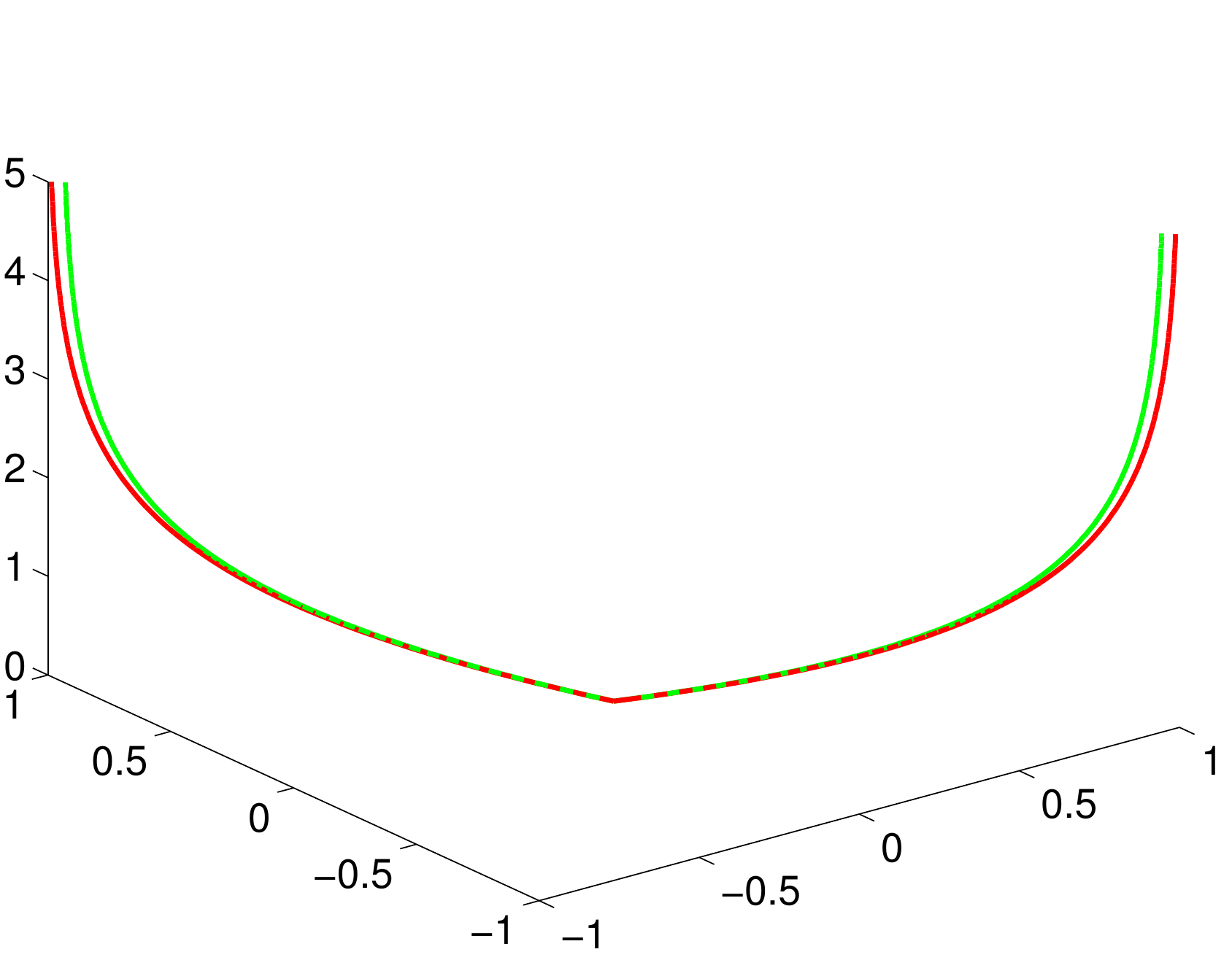}
\end{center}
\caption{reachable sets and minimum time functions for Example~\ref{ex:counter_ex_1}}
\label{fig:counter_exam}
\end{figure}\\
The two red points are the end points of the line segment for a smaller time $t_1 = 1$,
the two blue points are the end points of the line segment for a larger time $t_2 = 2 > t_1$.
The blue line segment is the reachable set for time $t_2$ (also the reachable set up to
time $t_2$). \\
All four points are on the boundary of the blue set $\mathcal{R}(t_2)$, but the minimum time to reach
the two blue points is $t_2$, while the minimum time to reach
the two red points is $t_1 < t_2$ which is a contradiction to 
\cite[Proposition~2.19]{BLp1}. 
\end{Example} 

\subsection{Problematic examples}

The first two examples show linear systems with hidden stability properties so that the discrete
reachable sets converge to a bounded convex set if the time goes to infinity (or is large enough
in numerical experiments). For larger end times the numerical calculation gets more demanding,
since the step size must be chosen small enough according to~\cite[Proposition~3.10]{BLp1}.
The remaining part of the subsection demonstrates that a target or a control set not containing the origin 
 (as a relative interior point) might lead to non-monotone behavior
of the (union of) reachable sets. In all of these examples the union of reachable sets is no longer convex.

\begin{Example} 
\label{ex:n1}
We consider the following time-dependent linear dynamics:
\begin{equation}\label{examplen1}  
\dot{x}_1=-x_2,\,\,\dot{x}_2=x_1 - \frac{1}{t^2} u,\,\,u\in [-1,1]
\end{equation}
The reachable sets converge towards a final, bounded, convex set due to
the scaling factor $\frac{1}{t^2}$ in the matrix $B(t)$, see~Fig.~\ref{fig:exam_n1}~(left).
From a formal point of view the strict expanding condition~\cite[(3.25) in Proposition~3.10]{BLp1}
is satisfied, but the positive number $\varepsilon$ tends to zero for increasing end time. On the other hand
we would stop the calculations if the Hausdorff distance of two consequent discrete reachable sets 
is below a certain threshold.
\begin{figure}[htp]
\begin{center}
\includegraphics[scale=0.4]{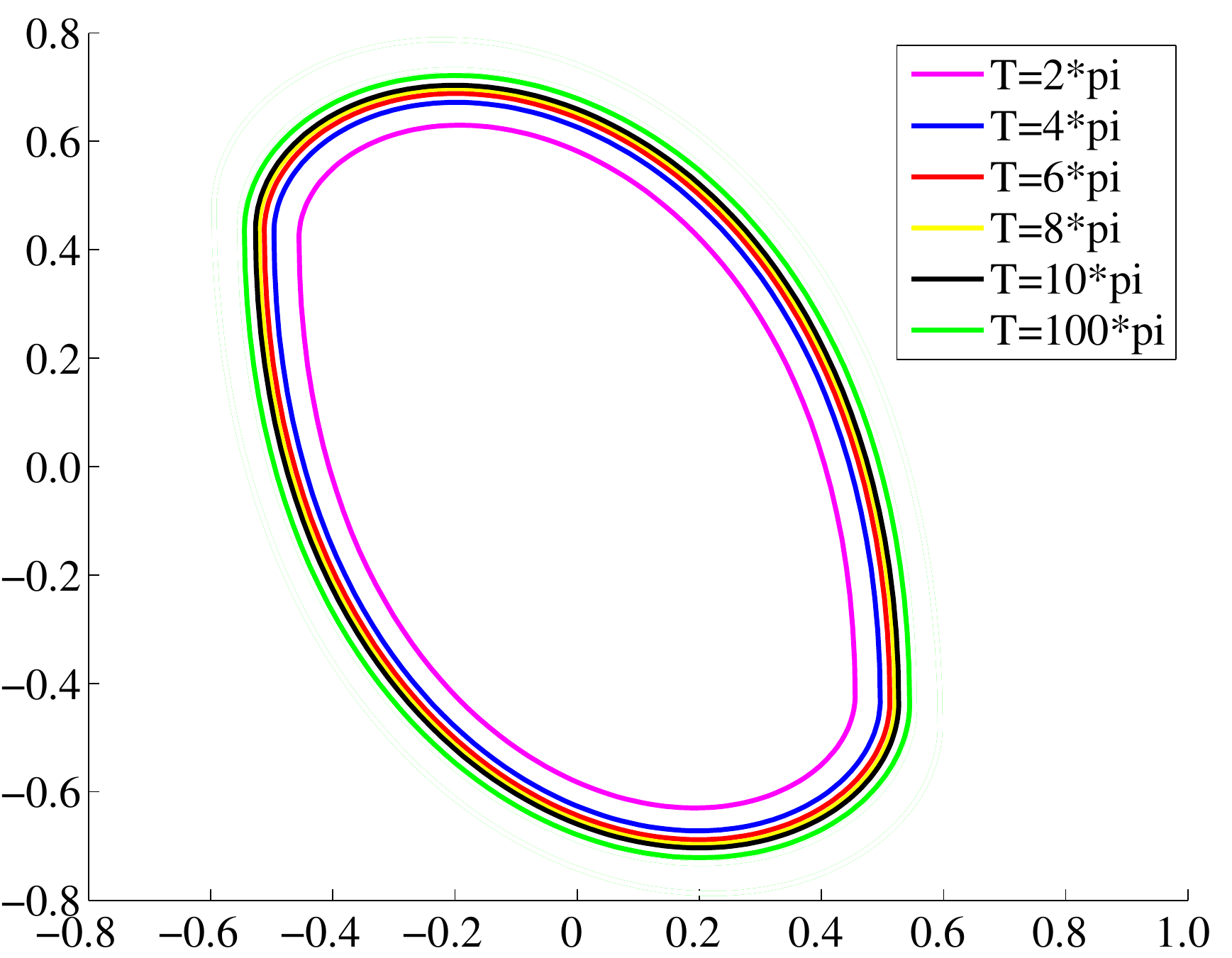}
\ \,
\includegraphics[scale=0.4]{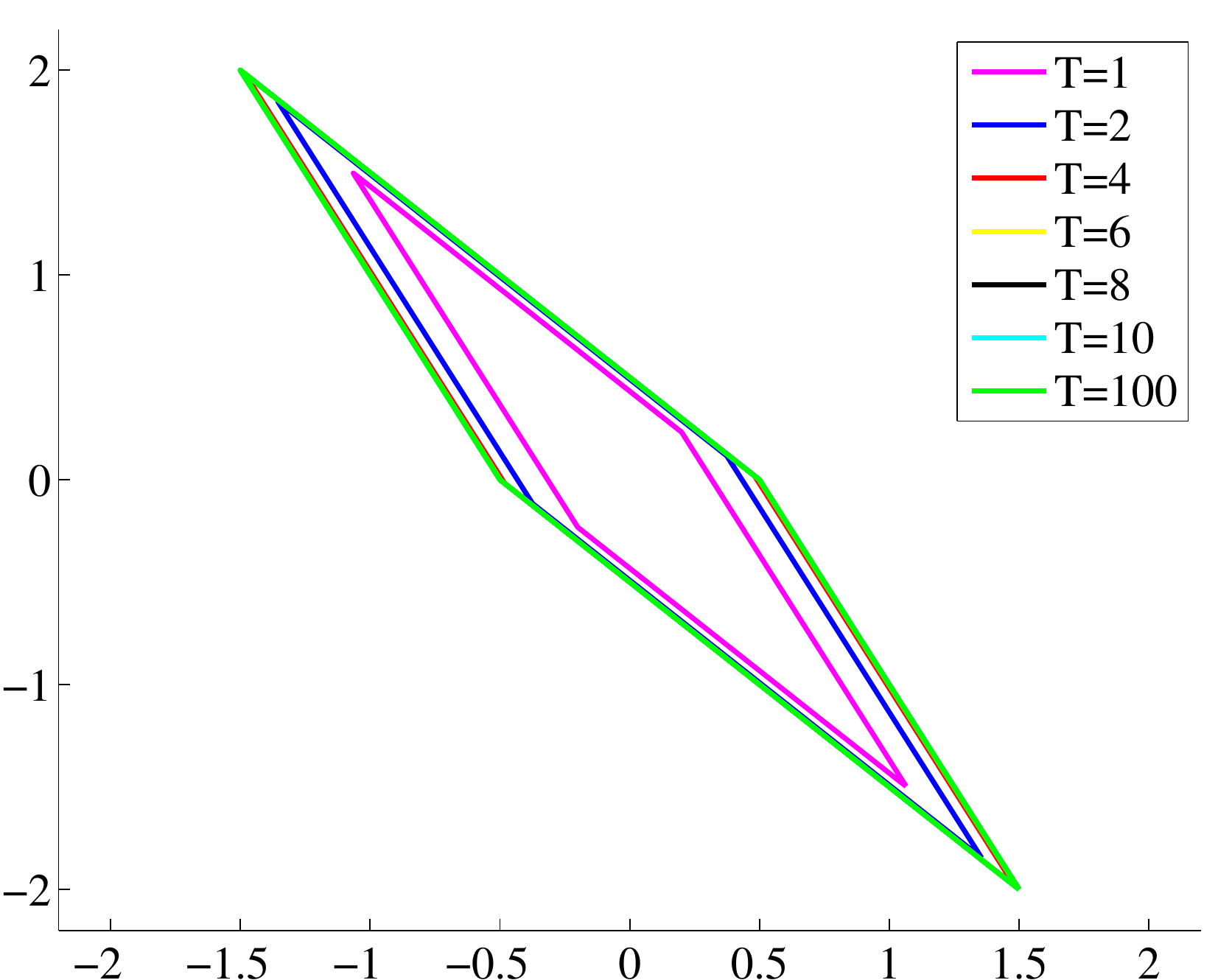}
\caption{reachable sets with various end times $t_f$ for Examples~\ref{ex:n1} and \ref{ex:n2}} 
\label{fig:exam_n1}
\end{center}
\end{figure}
\end{Example} 

\begin{Example}\label{ex:n2} 
   We reconsider Example~\ref{ex:3a} on the larger time interval $[t_0,t_f]=[0,100]$. $\bar{A}$ has negative eigenvalues -1 and -2. 
   Hence, the reachable sets converge towards a final, bounded, convex set, 
   see~Fig.~\ref{fig:exam_n1}~(right). We experience the same numerical problems as 
   in Example~\ref{ex:n1}.
\end{Example} 

\begin{Example} 
\label{ex:n4}
   Let the dynamics be given by
\begin{equation}\label{examplen4}  
\dot{x}_1=x_2 + u_1,\,\,\dot{x}_2=-x_1 + u_2,\,\,u\in B_1(0).
\end{equation}
In case a) the reachable sets for a given end time are always balls around the origin
(see~Fig.~\ref{fig:exam_n4} (left)),
if the target set is chosen as the origin.
In case b) the point $(2,2)^\top$ is considered as target set. Fig.~\ref{fig:exam_n4}~(right) shows
that the union of reachable sets is no longer convex.
\begin{figure}[htp]
\begin{center}
\includegraphics[scale=0.45]{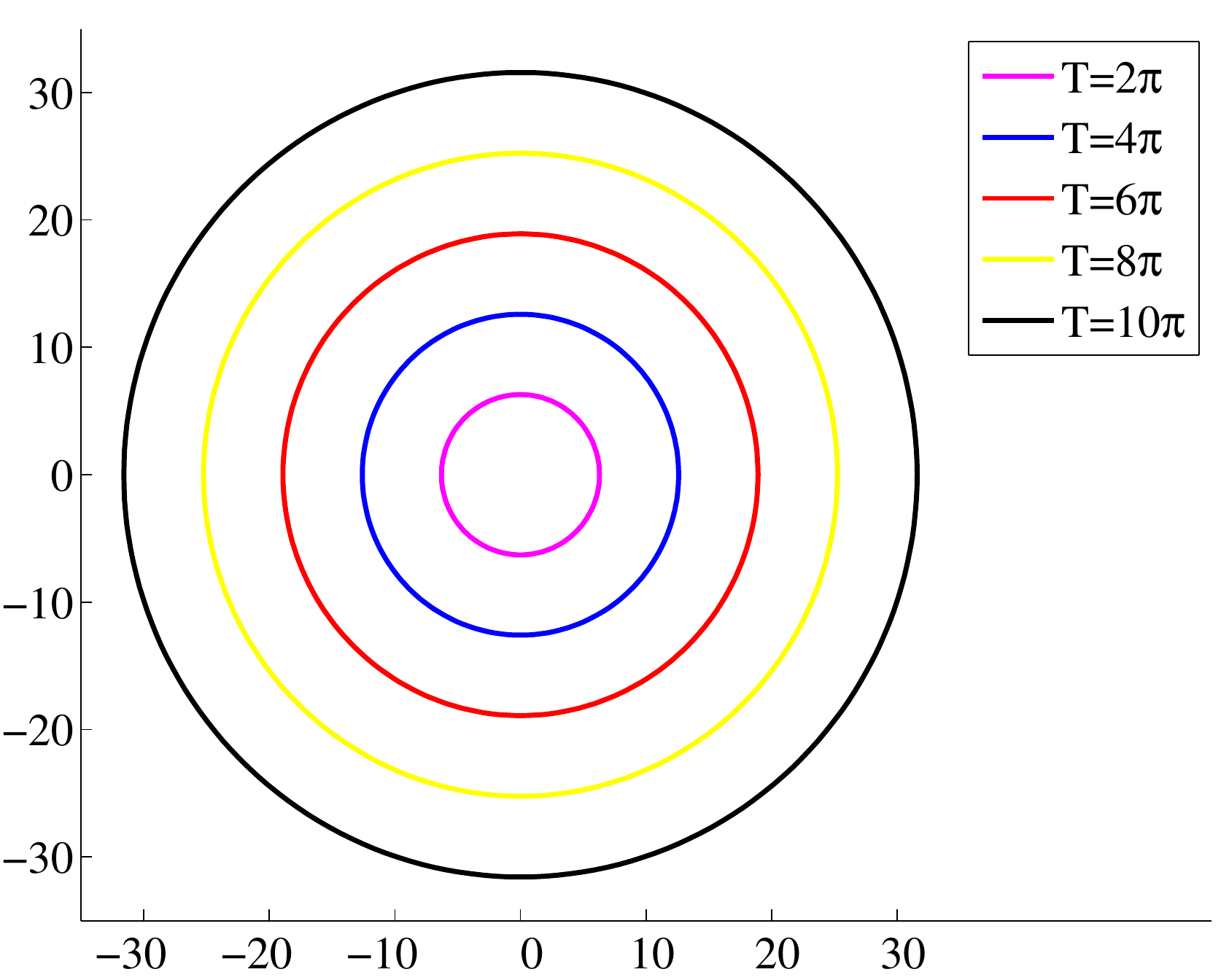}
\quad
\includegraphics[scale=0.45]{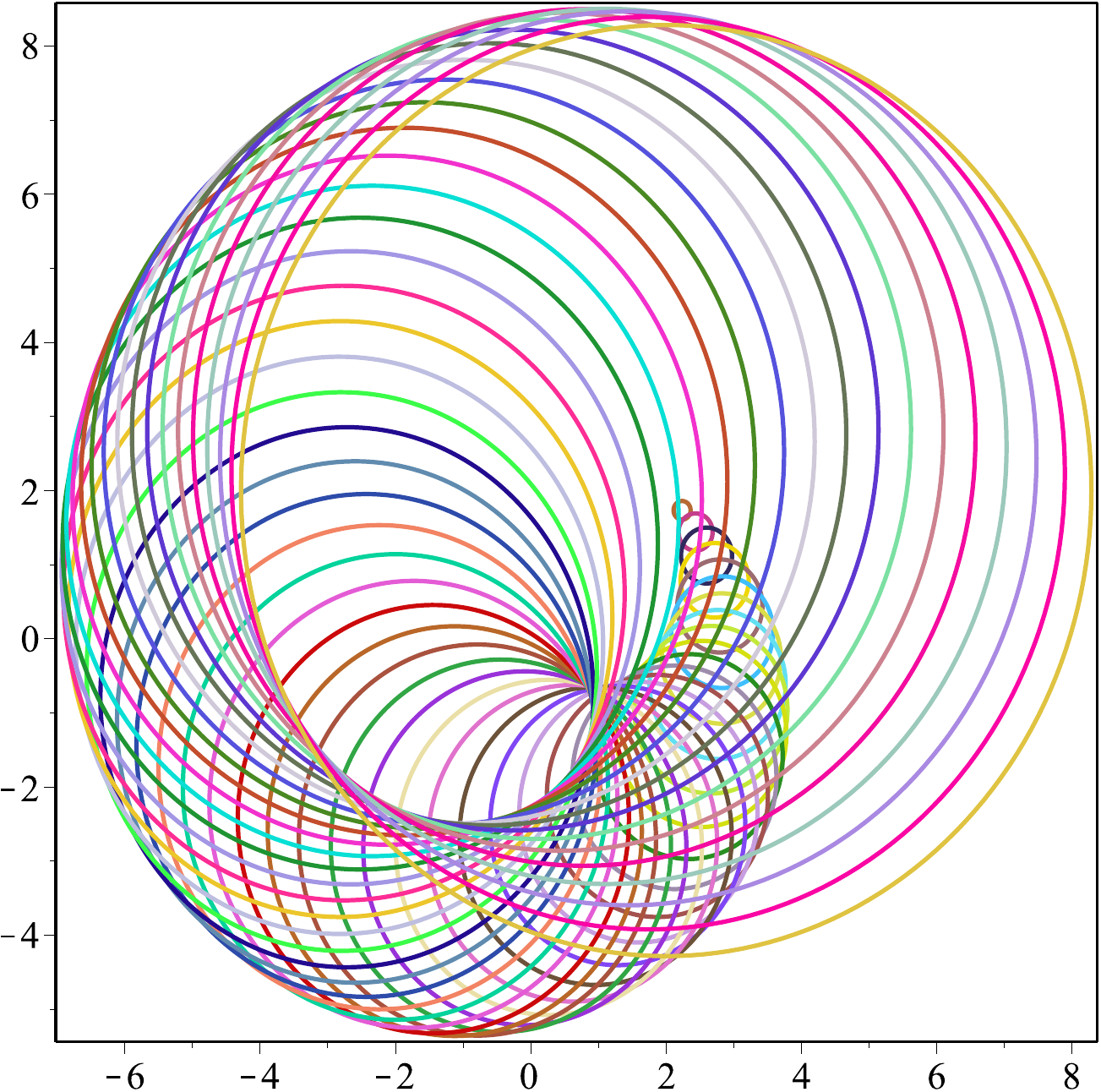}
\caption{reachable sets with various end times and different target sets for Example~\ref{ex:n4}} 
\label{fig:exam_n4}
\end{center}
\end{figure}
\end{Example} 

\begin{Example}\label{ex:n5}
Let us reconsider the dynamics \eqref{example3} of Example \ref{ex:2}, i.e.
\begin{equation*}
\dot{x}_1=x_2,\,\dot{x}_2=u,\,\,u\in U.
\end{equation*}
In the first case, let $\mathcal{S}=\bb{0},\,\,U=[0,1]$. From the numerical calculations, we observe that $\mathcal{R}(t)$, $\mathcal{R}_{\le}(t)$ are still convex and satisfy \cite[(2.12) in Remark 2.20]{BLp1}, but violate the strictly expanding property \eqref{eq:strict_exp_prop} as shown in Fig.~\ref{fig:exam_n5}~(left). In the other case, $U=[1,2]$ is chosen. The convex reachable set $\mathcal{R}(t)$ is not only enlarging, but also moving which results in the nonconvexity of $\mathcal{R}_{\le}(t)$. Moreover, both \cite[(2.12) in Remark 2.20]{BLp1} and \eqref{eq:strict_exp_prop} are not fulfilled in this example as depicted in Fig.~\ref{fig:exam_n5}~(right).
\begin{figure}[htp]
\begin{center}
\includegraphics[scale=0.4]{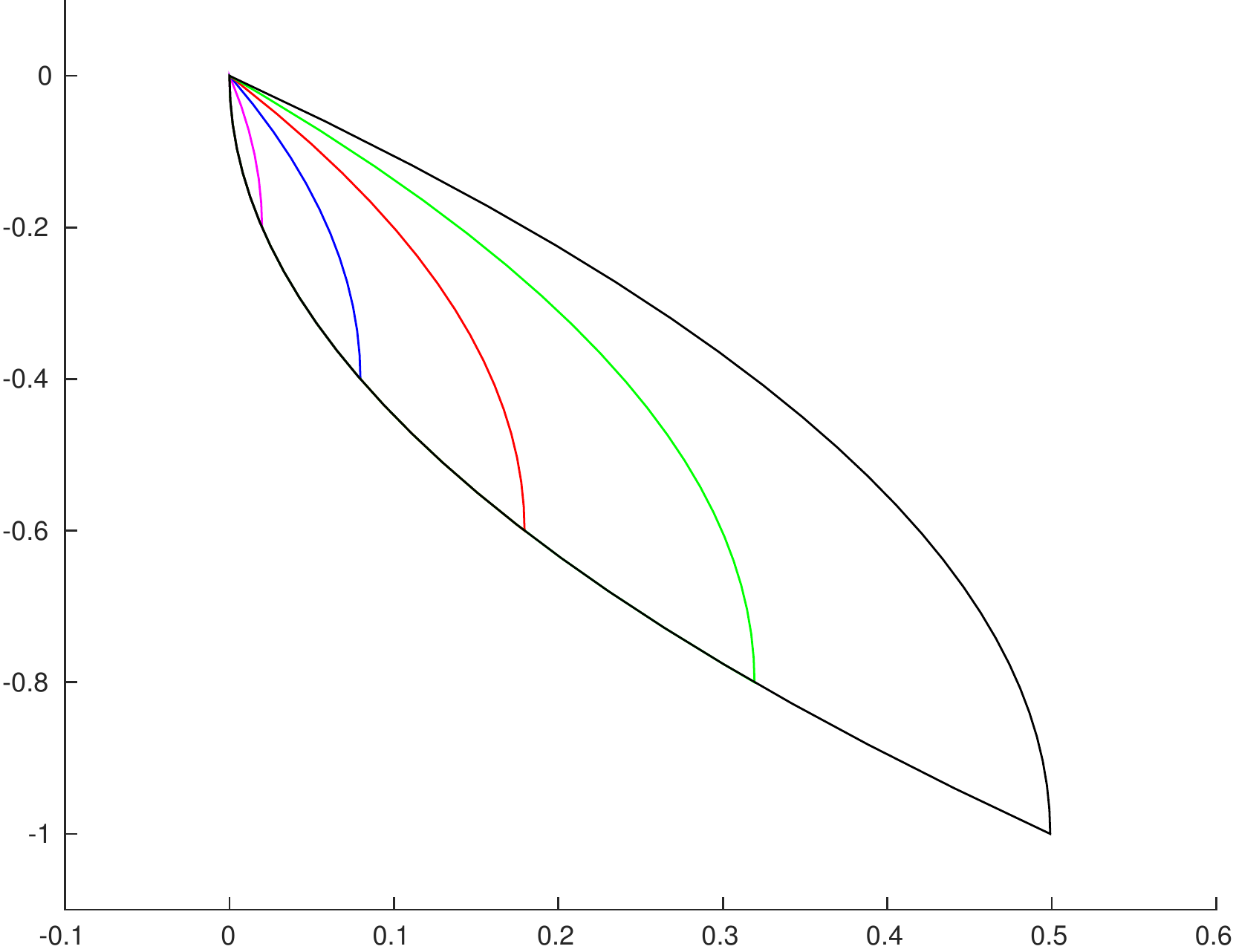}
\quad
\includegraphics[scale=0.4]{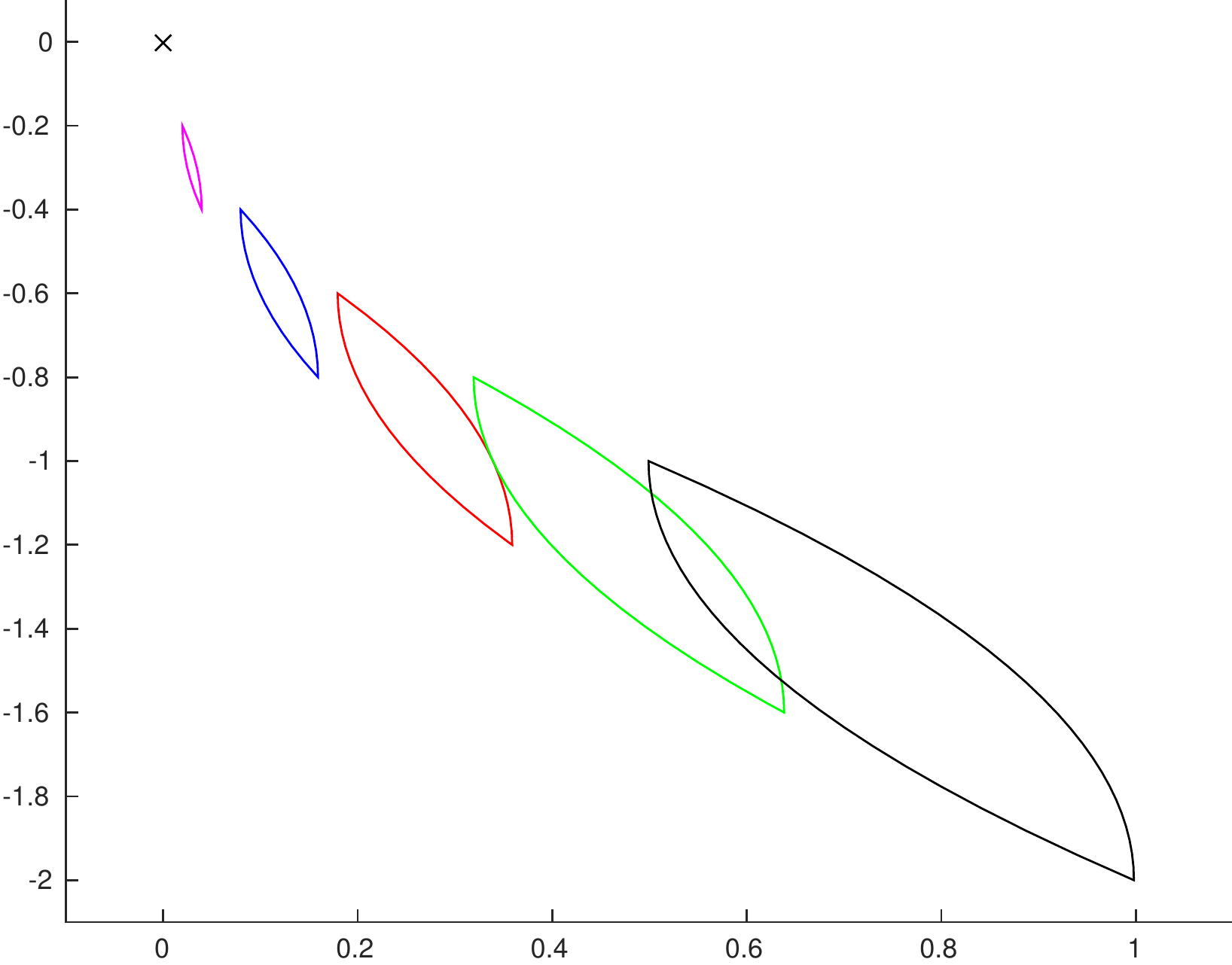}
\caption{reachable sets with various end times and different control sets for Example~\ref{ex:n5}} 
\label{fig:exam_n5}
\end{center}
\end{figure}
\end{Example}




\section{Conclusions}
\label{sec:concl}

Although the underlying set-valued method approximating reachable sets in linear control 
problems is very efficient,
the numerical implementation is a first realization only and can still be considerably improved.
Especially, step~3 in~\cite[Algorithm~3.4]{BLp1} 
can be computed more efficiently as in our
test implementation. Furthermore, higher order methods like the set-valued Simpson's rule
combined with the Runge-Kutta(4) method are an interesting option in examples where the underlying
reachable sets can be computed with higher order of convergence than 2, especially if 
the minimum time function is Lipschitz. But even if it is merely H\"older-continuous
with $\frac{1}{2}$, the higher order in the set-valued quadrature method
can balance the missing regularity of the minimum time function and improves the error
estimate.
We are currently working on extending this first approach for linear control problems
without the monotonicity assumption on reachable sets and for nonlinear control problems.

\section*{Acknowledgements}
The authors want to express their thanks to Giovanni Colombo and Lars Gr\"une
who supported us with helpful suggestions and motivating questions.


\end{document}